\documentclass[a4paper,final]{article}
\usepackage{amsmath}
\usepackage[numbers]{natbib}
\usepackage{amsthm}
\usepackage{amssymb}
\usepackage{comment}
\usepackage{graphicx}
\usepackage{multirow}
\usepackage{longtable}
\usepackage{wrapfig}
\usepackage{float}
\usepackage{xcolor}
\usepackage{hyperref}
\usepackage{setspace}
\usepackage{mathtools}
\usepackage{algorithm}
\usepackage{algorithmic}
\renewcommand{\algorithmiccomment}[1]{\bgroup\hfill//~#1\egroup}
\usepackage{eurosym}
\usepackage[utf8]{inputenc}
\usepackage[english]{babel}
\usepackage{amsfonts}
\usepackage{authblk}
\usepackage{xcolor}
\usepackage[normalem]{ulem}

\title{A decision support system for optimised industrial water management}

\author[1,*]{Stavros Vatikiotis}
\author[1]{Ioannis Avgerinos}
\author[2]{Stathis Plitsos}
\author[1]{Georgios Zois}
\affil[1]{\small ELTRUN Lab, Department of Management Science and Technology, Athens University of Economics and Business,\protect\\ \small Patision 76, 10434, Athens, Greece }

\affil[2]{\small University of Piraeus, Department of Industrial Management \& Technology,\protect\\ \small 80 Karaoli \& Dimitriou str, 18534, Piraeus, Greece 

\\ E-mail: {\tt stvatikiotis@aueb.gr, iavgerinos@aueb.gr, stathisp@unipi.gr, georzois@aueb.gr}}

\affil[*]{Corresponding Author} 
\date{} % clear date

\begin{document}
\maketitle

\begin{abstract}
Water scarcity and the low quality of wastewater produced in industrial applications present significant challenges, particularly in managing fresh water intake and reusing residual quantities. These issues affect various industries, compelling plant owners and managers to optimise water resources within their process networks. To address this cross-sector business requirement, we propose a Decision Support System (DSS) designed to capture key network components, such as inlet streams, processes, and outlet streams. Data provided to the DSS are exploited by an optimisation module, which supports both network design and operational decisions. This module is coupled with a generic mixed-integer nonlinear programming (MINLP) model, which is linearised into a compact mixed-integer linear programming (MILP) formulation capable of delivering fast optimal solutions across various network designs and input parameterisations. Additionally, a Constraint Programming (CP) approach is incorporated to handle nonlinear expressions through straightforward modeling. This state-of-the-art generalised framework enables broad applicability across a wide range of real-world scenarios, setting it apart from the conventional reliance on customised solutions designed for specific use cases. The proposed framework was tested on 500 synthetic data instances inspired by historical data from three case studies. The obtained results confirm the validity, computational competence and practical impact of our approach both among their operational and network design phases, demonstrating significant improvements over current practices. Notably, the proposed approach achieved a 17.6\% reduction in freshwater intake in a chemical industry case and facilitated the reuse of nearly 90\% of wastewater in an oil refinery case.
\end{abstract}


\section{Introduction} \label{Section:Problem}
Water scarcity is a pressing issue, affecting approximately 40\% of the global population due to insufficient water quantity, especially when considering water quality as well \cite{van2021global}. The root causes stem from inadequate sanitation and waste management services in many urban and rural areas, leading to significant pollution of water bodies \cite{spira2021water}. Factors contributing to this crisis include population growth, rising food demands, urbanisation, climate change \cite{levinson2008global}, and economic development \cite{figueiredo2021water}. Driven by these factors, water use has been annually increasing by 1 \% over the last 40 years, while it is expected to continue on a similar rate up to 2050 \cite{unesco2023}. Addressing this challenge requires a comprehensive approach that identifies consumption patterns and loss points, especially considering that only a tiny fraction of Earth's water is suitable for human use \cite{jha2018smart}. Industrial production plays a significant role in this matter. According to \cite{fao2022}, it is estimated that 45\% of the total abstracted water in Europe is utilised on industrial applications. Consequently, the industrial sector bears a substantial responsibility to tackle these challenges through technological advancements \cite{sengupta2017industrial}.

The above estimations, when combined with severe water scarcity problems on multiple areas of the world, encourage initiatives to efficiently deal with both situations, aiming at reducing the fresh water quantities drawn. As a result, the problem faced is two-fold: a) to efficiently manage fresh water resources without hindering the production processes inside the plant and b) to treat wastewater flows with the goal of reusing them, and thus, indirectly contribute to the first goal of fresh water intake minimisation. The complexity of the above problem calls for the design of an \textit{intelligent module} capable to decide upon optimal water flows within a process network on both operational and design time. That module has to be encapsuled within a \textit{decision support system (DSS)} in order to facilitate the work and choices to be made by the appropriate decision makers on the plant.

\noindent\textbf{Contribution.} 
%Our contribution comes forward on multiple ways. Firstly, we provide the user requirements of a service that assists in (waste-)water management optimisation and simulation. This is accompanied with a generic data model representing any (or the at least the majority) network of similar characteristics structure along with the respective JSON structure. Secondly, we provide an elegant and generic mathematical programming formulation compared to the ones presented in the literature, where the linearisation proposed by \cite{Lodi15} is used, in order to avoid solving non-linear formulations. As we will show on the next sections, the limitations of this linearisation are dealt with, by incorporating equivalent network representations (dummy nodes). Finally, we validate the efficiency of our methodology by applying the service as a whole on three case studies from different industries and with different scopes from an EU-funded project.
Our contribution emerges in multiple ways. We address the user requirements of a service that assists in industrial water management. This is accompanied with a generic data model representing networks of similar characteristics. Moreover, a quite generic and compact Mixed Integer Nonlinear Programming (MINLP) formulation is proposed, contrary to the ones presented in the literature (see Section \ref{sec:literature}), together with a linearisation procedure that efficiently tackles nonlinear constraints, inspired by \cite{Lodi15}. Interestingly, as shown on Section \ref{Section:Networks}, we overcome the limitations of the latter linearisation approach \cite{Lodi15} via a reduction to equivalent network representations (dummy nodes). The coupling of the proposed data model and MILP formulation allows us to deal with different objectives, such as cost, flow, and energy consumption, and different decision support phases, i.e., the design and operational phase. As a result, our versatile approach enables an industrial user to both examine different network scenarios on the design phase, and then, employ the corresponding decisions on the operational phase for different optimisation criteria. Also, a Constraint Programming (CP) formulation, which is equivalent with the MINLP model, is presented; nevertheless, its comparison with the linearised model showed inferior performance. Finally, the efficiency of our methodology is validated by applying the service as a whole on three case studies from different industrial fields, including an oil refinery and two chemical industries.% from an EU-funded project.

\noindent\textbf{Outline.} The remainder of this paper is structured as follows. In Section \ref{sec:literature}, an extensive literature review is presented. In Section \ref{Section:CaseStudies}, we provide a detailed description of the case studies. In Section \ref{Section:Networks}, we describe the optimisation method proposed. In Section \ref{Section:DSSAnalysis}, the analysis of the DSS is provided including user and data requirements with the system architecture. In Section \ref{Section:DSSUseScenaria}, we describe the scenarios under which the DSS is used, and experimentally demonstrate the benefits of applying the DSS. Conclusions are provided in Section \ref{Section:conlcusion}.
\section{Literature review} \label{sec:literature}
Water serves various purposes in the industrial sector, playing a role in numerous processes. It is employed for tasks such as product preparation and usage, cleaning, machinery cooling and heating, employee use \cite{rao2015energy}, and the treatment of raw materials \cite{walsh2016industrial}. We split the literature review in two parts: a) first we focus on - the commonly overlooked matter of - industrial water management DSSs, and b) we provide a detailed review specifically on optimisation methods that were applied in such a context.

\subsection{Literature Review on water management DSSs}

Regarding research on the DSS front, it is quite interesting that, across the literature most approaches focus on urban water management \cite{jha2018smart, owen2018smart, wang2015cyber}. Nevertheless, there is a considerable effort put on industrial water management approaches, encompassing systems, tools, or frameworks, that have been devised and implemented in industries, each tailored to specific needs and capabilities.

In \cite{atoui2021coupling}, the authors concentrate on the water heating process, employing a Bayesian classifier that combines statistical decisions and a fault signature matrix. This approach allows them to detect faults in the process by monitoring water flow and temperature. In \cite{ji2021research}, a customised water management system designed for a thermal power plant is presented. This system, equipped with monitoring and control capabilities, facilitates adaptive optimisation of operating parameters affecting water resource utilisation and the stability of the plant's water system. Finally, \cite{gupta2020smart} emphasise on the significance of monitoring and reporting capabilities in industrial water management systems. The authors stress the importance of addressing issues such as water losses, maintenance costs for the water network (including pipelines and valves damaged by water bursts) and electricity costs associated with water management.

\subsection{Literature Review on water management optimisation methods}

Regarding the optimisation front, a vast number of research works are focused on the modeling of the processes inside the process networks \cite{Dimartino21, Gu23, Kim20, Zhang19}. However, since the problem presented here considers the network on a system's level and not on the interior of each process, we focus only on relative research where the underlying problems involve the optimisation from the network flow perspective. 

We begin our literature review by examining works focused solely on wastewater inlet streams. \cite{Ang18} study various wastewater treatment network structures (series, parallel, and mixed), minimising costs, and find parallel processes most effective based on random data. \cite{Izquierdo08} design networks considering pipe diameters and slopes, using Particle Swarm Optimisation (PSO) to minimise costs despite energy constraint nonlinearities. PSO is applied to 60 network instances, achieving solutions for 59. \cite{Puchongkawarin15} model a wastewater network as a Mixed Integer Nonlinear Program (MINLP) to assess the economic sustainability of wastewater reuse in a wine distillery case study. On the other hand, some applications require the consideration of uncertainty parameters. As a result, \cite{Savun-Hekimoglu23} develop a MILP for network design in an industrial park, focusing on Net Present Value (NPV) under uncertain discharge amounts over a 30-year horizon. Similarly, \cite{Tosarkani20} propose a robust flexible chance-constrained model (RFCCM) for identifying cost-efficient network designs considering uncertain treatment process costs or available wastewater, applied to shale gas production. \cite{Alfaisal23} minimise sewer network design costs by selecting component locations under capacity uncertainty using a MILP, tested on hypothetical data. Additionally, several studies address multi-objective optimisation in their problem formulations. \cite{Hreiz15} optimise operating costs and nitrogen discharge through a Genetic Algorithm (GA), exploring excessive sludge reuse for energy production. \cite{Niu22} optimise energy consumption and effluent quality using deep learning and dynamic multi-objective ant lion optimisation (DMOALO) for solving a benchmark instance. \cite{Perez2023} add social aspects (i.e. job offerings) and address wastewater quantity uncertainties in a MILP, tested on a case study with high computational times but offering practical insights. However, these complex, case-specific formulations limit transferability. In contrast, \cite{Bozkurt15} propose a relatively more generic MINLP for early-stage treatment network design, tested on 50 random instances with three cost-related objectives.

In addition to wastewater treatment network design, significant research also focuses on managing and optimising freshwater inlet streams. \cite{Abdulbaki17} address the allocation of freshwater resources to meet varying demands across multiple industrial plants, aiming to minimise costs. They propose a Mixed-Integer Nonlinear Programming (MINLP) model, which is linearised and applied to a small city case study under five scenarios to demonstrate the model's practical application. Further studies have tackled problems involving different units, such as processing and regeneration units, each with different objectives like minimising freshwater intake or optimising interconnections between system components. Notable works to this direction include \cite{Boix11} and \cite{Ramos14}, who use MINLP models for similar goals in industrial settings. Finally, \cite{Cassiolato21} extend this research with a nonlinear model to minimise construction costs of water distribution networks, applying it to four case studies from the literature and yielding improved results.

However, in industrial applications, it is often necessary to consider both freshwater and wastewater sources simultaneously. \cite{Ortega09} explore this concept by mixing different inlet streams, which are then directed to a set of processes designed either for recovery or discharge, while ensuring compliance with strict quality limits. Similarly, \cite{Yang14} and \cite{Yang15} apply a Mixed-Integer Nonlinear Programming (MINLP) model, coupled with a decomposition formulation, to minimise the total treatment cost across the network. \cite{Pungthong15} tackle a similar problem using a two-stage NLP-based heuristic, first minimising freshwater intake and subsequently reducing the total treatment cost. \cite{Al-Zahrani16} focus on transferring multiple water resources to various applications (e.g., industry, agriculture) with objectives related to both cost and resource quantities, such as wastewater reuse and groundwater minimisation. This approach is applied in Saudi Arabia on 100 synthetic data instances using a Goal Programming (GP) model. In another context, \cite{Hansen18} develop a mathematical programming formulation for a petrochemical plant to efficiently treat wastewater while minimising freshwater intake and meeting the required flowrate demands. \cite{Gaurav24} review three case studies from the literature, proposing a nonlinear program (NLP) that significantly reduces freshwater consumption. Finally, \cite{Cao23} introduce a Mixed-Integer Linear Programming (MILP) model tailored for shale gas production, showing substantial benefits in wastewater treatment and freshwater minimisation.

\subsection{Literature Review Summary}

%However, it should be noted that the majority of the aforementioned research works rely on nonlinear solvers, which hinder the scalability of the application of the proposed models. Furthermore, the designed mathematical formulations are case-specific, tailored on the proposed case studies, and thus, lacking the capability of generalising on different cases or applications.

To sum up, Table \ref{tab:review} provides a comparison between the research works of closely-related literature. In order to simplify the notation, works aiming at cost minimisation are denoted by $TC$, $C$ is related to works where the total number of connections within the network is minimised and $FW$, $WW$ represents the minimisation of fresh water intake and maximisation of treated wastewater respectively. Additionally, $Q$ denotes cases where quality parameter values are included in the objective. `General Model' refers to literature works, where the proposed model could be adopted to other potential cases. Approaches used on the design phase are denoted as $D$, while the ones examining the problem on the operational phase are denoted by $O$. Optimisation methodology is presented on columns \textit{Exact Opt. Method} and $(Meta-)heuristics$. Column $DSS$ shows whether a DSS approach has been provided in the corresponding research paper. Finally, column $\# Cases$ denotes the number of case studies examined. Note that, well known benchmark models or previously reported case studies are not considered on this column.

\begin{table}[H]
\centering
\scriptsize
\caption{Literature Review Table}
\resizebox{\textwidth}{!}{
\begin{tabular}{cccccccccc}
\hline
Reference        & Objective                  & Inlet  & Generic Opt.   & Phase  & Opt. & $DSS$ & \# Cases\\
        &                   &   &  Model  &   & Method  &  & \\ \hline

\cite{Izquierdo08} & $TC$                  & WW        & -                    & D & PSO    & -      &  - \\
\cite{Ortega09}    & $ TC $  & WW, FW               & -             & D & MINLP                   & -      & - \\
\cite{Boix11}    & $ FW, C $  & FW               & -             & D & MINLP          &     -  & -    \\ 
\cite{Ramos14}    & $ FW, C $  & FW               & -             & D & MINLP           &      - & -  \\ 
\cite{Yang14}    & $ TC $  & WW, FW               & -             & D & DC               & -          & 1 \\
\cite{Pungthong15}    & $ TC $  & WW, FW               & -             & D & MINLP            & -           & 1 \\
\cite{Puchongkawarin15}    & $ TC $  & WW  & -     & D & MINLP            & - &     1    \\ 
\cite{Yang15}    & $ TC $  & WW, FW               & -             & D & MILP         & -             & - \\
\cite{Hreiz15}    & $ TC, Q $  & WW  & -     & D &   GA   &     -  & -   \\ 
\cite{Bozkurt15}    & $ TC $  & WW               & \checkmark             & D & MINLP          & - & - \\
\cite{Al-Zahrani16}    & $ TC,FW,WW $  & WW,FW               & -             & D & GP          & -               & 1 \\
\cite{Abdulbaki17}    & $ TC $  & FW               & -             & D & MINLP, MILP                 & -      & 1 \\
\cite{Hansen18}    & $ FW $  & WW, FW               & -             & O & NLP               & -        & 1 \\
\cite{Ang18}    & $ TC $  & WW               & -             & D & MINLP              & -        & - \\
\cite{Tosarkani20}    & $ TC $  & WW              & -             & D & RO             & -           &  1\\
\cite{atoui2021coupling}    & -  & WW, FW              & -             & O & -           & \checkmark         & 1 \\
\cite{Cassiolato21}    & $ TC $  & FW              & -             & D & MINLP           & -              &  -\\
\cite{ji2021research} & -  & WW, FW              & -             & O & -         & \checkmark         & 1 \\
\cite{Niu22}    & $ EC, Q $  & WW              & -             & O    & DMOALO      & -         & - \\
\cite{Alfaisal23} & TC & WW & - & D & MILP& - & - \\
\cite{Cao23} & TC & WW,FW & - & D & MILP & - & 1\\
\cite{Perez2023} & TC, Q & WW & - & D & MILP &  - & 1 \\
\cite{Savun-Hekimoglu23} & TC & WW & - & D & MILP &  - & 1\\
\cite{Gaurav24} & TC & WW,FW & - & D & MINLP & - & - \\
Current Work    & $ TC,FW,WW $  & WW,FW               & \checkmark            & D, O & MINLP, MILP, CP    & \checkmark      &  3 \\

\hline
\end{tabular}
}
{\tiny NLP: Nonlinear Programming, MINLP: Mixed Integer Nonlinear Programming, DC: Decomposition, RO: Robust Optimisation, PSO: Particle Swarm optimisation, DMOALO: Dynamic Multi-objective Ant Lion Optimisation, GA: Genetic Algorithm}
\label{tab:review}
\end{table}

%we describe the wastewater management problem (WWMP) and design optimisation formulations which can be easily adopted to be applicable to several wastewater flow networks. Additionally, we establish possible applications of the presented formulations. Secondly, we present a high-level framework, where optimisation algorithms are interacting in real-time with different modules as data analytics, simulation and real time monitoring platforms. 

Table \ref{tab:review} outlines the key strengths of our contribution to the existing literature, providing a basis to emphasise several important observations. First, our work is unique in combining a DSS with an optimisation method. On one hand, an optimisation model alone is insufficient if it lacks an appropriate support system to enable end-user application. On the other hand, a DSS focused solely on monitoring assets within a plant cannot provide optimal solutions. Thus, combining a DSS with optimisation is crucial for addressing industrial challenges. Moreover, the implementation of a generic optimisation model in our approach allows for a multifaceted examination of the problem, accommodating diverse objectives, constraints, and phases, e.g., operational and design. From a technical standpoint, it is important to note that most related works rely on nonlinear solvers, which limit the scalability of their models. In contrast, our mathematical formulation is linear and intentionally designed to be case-independent, providing enhanced flexibility and applicability across diverse scenarios. Lastly, while most existing literature typically examines a single case study, our framework is validated through three distinct case studies, showcasing its broad applicability and robustness. \newline

\section{Water management in the chemical and oil-refining industries} \label{Section:CaseStudies}
In this section we outline our case study approach (Section \ref{Section:CS_Approach}) and furnish the specifics of the case studies (Section \ref{Section:CS_Descriptions}) derived from three industrial organisations actively participating in an ongoing EU-funded research project dedicated to industrial water management. These organisations are motivated to engage in the project with the aim of exploring and implementing digital methods to reduce water consumption, facilitate water reuse, and mitigate water losses. In pursuit of these goals, they have collaborated with us, addressing challenges related to water management optimisation, the specific objectives they seek to improve, the respective constraints they have and the user requirements that facilitate the implementation of a service that optimises water management decisions in a plant.

\subsection{Case study approach} \label{Section:CS_Approach}
We have embraced a multiple case study approach to present diverse perspectives on the issue of water management in industrial sites. We opted for this methodology due to the presence of clearly defined cases with distinct boundaries, aiming to gain a comprehensive understanding of each case \cite{kegler2019study}. Among the three case studies, there are two chemical industries that share similarities in water usage but differ in terms of products, user requirements, and optimisation objectives. The remaining case study involves an oil refinery characterised by different water usage patterns, while having the same optimisation objective with one of the chemical industries. 

In all case studies the quality of the water used is significantly impacted, necessitating a corresponding water treatment process. However, there are different water network structures and divergent approaches to reusing or disposing the treated water. Consequently, our collection of case studies offers complementary insights in terms of water use and treatment, objectives, constraints and user requirements. This diversity enables us to generalise our findings across various industrial contexts. 

Semi-structured interviews were carried out for each case study, involving discussions with representative employees at each industrial site. These individuals are actively engaged with the existing applications employed for water efficiency objectives and possess awareness of the additional functionalities required to enhance the water-related objectives prioritised by the company. The interviews were conducted twice, involving 5–15 participants from the respective organizations, with each session lasting 2 hours. The roles and seniority levels of the employees involved are detailed for each case study in the respective sections below.

\subsection{Case study descriptions} \label{Section:CS_Descriptions}
In this section we describe how the aforementioned methodology is utilised on each of the case studies. 

\subsubsection{Oil refinery}
Oil \& gas refining industry is highly water intensive, requiring vast amounts of water, used as cooling water, service water, firefighting water, demineralisation water and water for steam production. Our case study serves as a representative example of this sector as it operates four oil refineries of a total annual processing capacity of 30 million tons crude oil. In this work, we focus on one of its refineries which consumes both fresh water from a lake and treated wastewater from its own wastewater recovery plants. The semi-structured interviews included participants from various roles, such as senior Process Engineers, Water Treatment Specialists, Environmental Engineers, Laboratory Technicians, and management personnel, including positions like Chief Information Security Officer and IT Infrastructure Coordinator.

The treated industrial wastewater is used as firefighting water and cooling water make up (with a flow rate of 350 $m^3/h$) and aims to increase the water reuse opportunities and decrease the fresh water intake from the lake, as any attempt approaching near zero discharge goal should be considered. There are two independent water tanks whose contents cannot be recycled or discharged to the receiving water body. The first one (T\_1) contains highly conductive backwash water (with a flowrate capacity of 300 $m^3/h$). The network includes also a second contaminated tank (T\_2) with the capability of mixing the two, in order to discharge the blended contaminated water to the receiving water body. The second tank (T\_2) contains hot wash water with a capacity of 110 $m^3/h$. The goal of this case study is to support network design decisions, including additional treatment components, as well as operational ones, maximising the treated waste water flowrate across a given time horizon.

\subsubsection{Chemical industry A}
This industry connects chemistry and innovation to the principles of sustainability. This involves specialty in chemicals, advanced materials, and plastics. The organisation offers technology-based products and services to clients in about 160 countries and in high growth sectors like food and specialty packaging, industrial and consumer packaging, health and hygiene, electronics, energy, architectural and industrial coatings, home and personal care, as well as infrastructure. The location of the plant in focus has a tight water balance, therefore characterised as water scarce. Hence, sustainable solutions for the responsible use of water resources are sought for long-term development. One such solution is the intake of fresh raw water streams to produce customised water qualities for different downstream applications (demineralisation of water and cooling tower make-up). The semi-structured interviews included participants from various roles, such as senior Water Specialists, Process Engineers and Environmental Engineers.

The identified goal is an implementation plan with the potential to further reduce the fresh water intake by an additional 20\%. All the above will move the plant towards the near-zero discharge goal, by decreasing the fresh water intake from the lake. Hence, the goal in this case study is the support on the network structure design that includes treatment units and additional water streams and the minimisation of fresh-water intake across a given time horizon for its operations.

\subsubsection{Chemical industry B} \label{subsection:chemb}
Chemical industry produces sodium carbonate, sodium bicarbonate (also for pharmaceutical use), calcium chloride, chlorine, hydrochloric acid, chloromethane, plastic materials, peracetic acid and hydrogen peroxide. The industrial wastewater resulting from the Peroxide and Peracetic Acid production amounts to about 10 $m^3$/h and is characterised by high concentrations of TOC/COD, nitrates, sulphates and, to a lesser extent, phosphates and hydrogen peroxide. Due to space limitations, which restrict the plant's capability of creating a complete inbound treatment plant, it has been decided a priori that a relatively smaller amount of wastewater will be treated inside the plant, while most of the purified water (needed for its operations) will be provided from a nearby water purification agent. The semi-structured interviews included participants from various roles, such as senior Process Engineers, Water Treatment Specialists and Environmental Engineers. The focus of this case study is placed on estimating the performance of a new network structure design, including this treatment process, and deciding on the optimal selection for water re-use by minimising the total cost.
\section{Problem definition and solution method} \label{Section:Networks}
Although the case studies described in Section \ref{Section:CaseStudies} appear divergent, they share a common core need for network optimization. Aiming at unifying diverse business requirements under the same framework, we refrain from designing separate formulations, resulting in the proposal of a holistic optimisation model. First, in Section \ref{Section:problemDefinition} the problem definition is provided, Section \ref{Section:mathematicalModel} presents a MINLP formulation of the problem, while in Section \ref{Section:linearApproximation}, a linearised MILP formulation is shown which accommodates the problem at hand.

%In this section we dive deeper into the method employed to optimise the operational and design decisions respectively. To do this, we provide in Section \ref{Section:problemDefinition} the problem definition, in Section \ref{Section:mathematicalModel} the mathematical modeling approach that is reflected in a MINLP formulation of the problem, and in Section \ref{Section:linearApproximation} a linearised MILP formulation that accommodates the problem at hand.

\subsection{Problem definition} \label{Section:problemDefinition}
Based on the above discussion, our optimisation problem deals with deciding upon the optimal flow rates across the system network. Before moving on to the mathematical modelling, we will first highlight the basic components included in a generic network which captures the main components and attributes of process networks formed within industries such as chemical and oil ones. %as it seems that a generic characterization of all nodes inside a network is possible. 
More precisely, three typical node categories are: a) \textbf{Flow provider nodes}, which provide flows to the network, but do not receive flows by other nodes, b) \textbf{Intermediate nodes}, which both receive and provide flows within the network, and c) \textbf{Flow receiver nodes},  which only receive flows, but do not provide any flows to other nodes. 

These categories encapsulate the different components of the network. Each component may be associated with specific attributes, inferred by the properties of each component, while inserted to the model based on user requirements and the availability of respective data. More precisely:

\begin{itemize}
    \item[-] \textbf{Fresh Water streams}, which insert fresh water to the network, accompanied by the attribute \textit{Quality parameters} (\textbf{$c_{ip}$}), the value of each quality parameter $p$ at the time the flow from $i$ is inserted in the network. 
    \item[-] \textbf{Wastewater streams}, which insert wastewater, accompanied by: i) \textit{Wastewater Flow} (\textbf{$q_i$}), a specific amount of already determined flow inserted by $i$ in the network, and ii) \textit{Quality parameters} (\textbf{$c_{ip}$}), the value of each quality parameter $p$ at the time the flow from $i$ is inserted in the network.
    \item[-] \textbf{Treatment processes}, which alter flow qualities, accompanied by: i) \textit{Flow rate reduction rate} (\textbf{$SR_i$}), the percentage of the (waste)water flow rate that exits each process compared to the inflow rate, ii) \textit{Flow rate fixed reduction} (\textbf{$SF_i$}), the fixed flow rate exiting $i$ regardless of the flow inserted, iii) \textit{Contamination reduction rate} (\textbf{$RR_{ip}$}), the reduction percentage of the value of quality parameter $p$ on each treatment process $i$, iv) \textit{Contamination fixed reduction} (\textbf{$RF_{ip}$}), the specified value of quality parameter $p$ when exiting $i$. 
        %\item[] \textit{Capacity} (\textbf{$C_i$}): The maximum amount of flow rate which can traverse via treatment process $i$. 
    \item[-] \textbf{Discharge points}, where (waste)water may be discharged, accompanied by \textit{Quality Requirements}, the minimum and maximum quality values between which water is eligible to be discharged.
    \item[-] \textbf{Applications}, referring to a set of end points (sinks), where water (or treated wastewater) may be used, accompanied by: i) \textit{Quality Requirements}, the minimum and maximum quality values between which water is eligible to enter each application, and 
ii) \textit{Demand}, the minimum water flow to direct to an application.
    \item[-] \textbf{Edges}, referring to the set of possible connections between different components inside the network, accompanied by \textit{Capacity}, the maximum flow rate which can traverse each edge.
\end{itemize}

Note that reduction rates ($RR$ and $SR$) receive values between 0 and 1, when the flow rate or the quality parameters are decreased, and greater than 1 if the aforementioned parameters are increased. As an example, if a flow rate $x$ is inserted in a component $i$ with a flow rate reduction rate $S_i$, then the output is: $y = S_i \cdot x$. Additionally, \textbf{Tanks}, i.e. components where (waste)water may be (temporarily) stored may be deemed as \textbf{Treatment Processes} with the values of $SR$ and $RR$ set to 1.

Apart from the aforementioned component specific attributes, we, also, include attributes which may be associated with more than one components: 
a) \textit{Capacity} (\textbf{$C_i$}), the maximum quantity to be drawn from $i$, b) \textit{Fixed Cost} (\textbf{$FC_i$}), the fixed cost when flows are drawn from or traverse $i$, c) \textit{Variable Cost} (\textbf{$VC_i$}), the variable cost (per $m^3/h$) for flows drawn from $i$, d) Fixed Energy (\textbf{$FC_i$}), the fixed energy consumption when $i$ is used, e) Variable Energy (\textbf{$VE_i$}), the variable energy consumption (per $m^3/h$). 

The next step is to associate each network component to an appropriate node category. Wastewater and fresh water streams are placed in the flow provider category. Similarly, discharge points and applications are considered as flow receiver nodes, while treatment processes and tanks are intermediate nodes.
However, it should be clarified that attributes may be extended to be associated with more components than the ones listed above. For example, if a treatment process is highly sensitive in terms of the inserted water quality, then, the quality requirements attribute may be also assigned to that process.

%Table \ref{tab:node_categories} sums up the relations between node categories and components. 

\begin{comment}

\begin{table}[ht]
\centering
\caption{Summary of node categories, network components and attributes}
\resizebox{\textwidth}{!}{
\begin{tabular}{ccccccc}
\hline
Node Category    & Components                 & Attributes\\ \hline

\multirow{3}{*}{Flow provider nodes} & Fresh Water Streams & Quality Parameters \\
& \multirow{2}{*}{Wastewater Streams} &      Flow          \\  
        & &  Quality Parameters \\
        \hline
\multirow{4}{*}{Intermediate Nodes} & Tanks &   Capacity     \\
        & \multirow{3}{*}{Treatment Processes}&  Flow Rate Reduction \\
        & &  Contamination reduction rate \\
         & &  Capacity \\
         \hline
\multirow{3}{*}{Flow receiver nodes} & \multirow{2}{*}{Applications} &   Quality Requirements             \\ 
        & &  Water flow demanded \\
         & Discharge &  Quality Requirements  \\
\hline
\end{tabular}}
\label{tab:node_categories}
\end{table}

\end{comment}

\begin{table}[ht!]
\scriptsize
    \centering
    {\resizebox{\textwidth}{!}{
    \begin{tabular}[t]{ll}
        \hline
        \textbf{Model Parameters} & \\
        \hline
        $I$             &   Set of Components\\
        $E$             &   Set of Edges\\
        $P$             &   Set of pollutants\\
        $d_{i}$         &   Flowrate demand of component $i\in I$ (in $m^3/h$)\\
        $q_{i}$         &   Fixed flowrate provided from component $i\in I$ to the network (in $m^3/h$)\\
        $C_{i}$         &   Maximum amount of flowrate to traverse component $i\in I$ (in $m^3/h$) \\
        $C_{e}$         &   Maximum amount of flowrate to traverse edge $e\in E$ (in $m^3/h$) \\
        $\delta^{i+}$   &   Subset of $E$; The edges in which component $i\in I$ is the origin node \\
        $\delta^{i-}$   &   Subset of $E$; The edges in which component $i\in I$ is the destination node \\
        $SR_{i}$         &   Flowrate reduction rate of component $i\in I$\\
        $SF_{i}$         &   Fixed Flowrate reduction of component $i\in I$\\
        $RR_{ip}$        &   Quality reduction rate in component $i\in I$ for pollutant $p\in P$ \\
        $RF_{ip}$        &   Fixed quality reduction in component $i\in I$ for pollutant $p\in P$ \\
        $l_{ip}$        &   Lower limit in component $i\in I$ for pollutant $p\in P$ \\
        $u_{ip}$        &   Upper limit in component $i\in I$ for pollutant $p\in P$ \\
        $\bar{c}_{ip}$  &   Given quality of pollutant $p\in P$ on a component $i\in I$ \\
        ${FC}_{i}$  &   Fixed Cost for the use of component $i \in I$\\
        $VC{i}$  &   Variable Cost per unit of flowrate traversing component $i \in I$\\
        $FE_{i}$  &  Fixed Energy consumption for the use of component $i \in I$\\
        $VE_{i}$  &    Variable Energy consumption per unit of flowrate traversing component $i \in I$ \\
        $M$  &  A sufficiently large number\\
        $\mu$  &  A sufficiently small number\\
        \hline
    \end{tabular}}}
    \caption{Model Parameters}
    \label{tab:annotation}
\end{table}

 Note also that the input data do not provide labels for each component. The pre-processing stage determines the labeling of $I$ as follows: i) If $\delta^{i-} = \text{\O}$, then component $i$ is a flow provider (a fixed flowrate $q_{i}$ is optional), ii) if $\delta^{i-}\neq \text{\O}$ and $\delta^{i+}\neq \text{\O}$, component $i$ is an intermediate node, iii) if $\delta^{i+} = \text{\O}$, then $i$ is a flow receiver, iv) if $|\delta^{i-}| = 1$, then $i$ receives flow from exactly one node.

Table \ref{tab:annotation} summarises the above discussion. Note that attribute $SR_i$ is in conflict with $SF_i$ on the same component $i \in I$, as there must be a single way to determine the outlet of each component. For example, a treatment process $i_1$ may inflict 20\% of water losses (thus, $SR_{i_1} = 0.8$), while a treatment process $i_2$ may always have an output of e.g. 300 $m^3/h$ (thus, $SF_{i_2} = 300$) regardless the entering flowrate. Similar conflicts are applied on attributes $RF_{ip}$ and $RR_{ip}$ per component and pollutant. For example, on a process $i_3$, the output of a pollutant $p_1$ may be equal to a value $c$ (thus, $RF_{i_{3}p_{1}} = c$), while the value of a pollutant $p_2$ is decreased by 40\% ($RR_{i_{3}p_{2}} = 0.6$). 
%{\color{red}+++ do we have to say sth more about the objectives of the problem - Answer: There is some discussion after the constraints.}

\subsection{Mathematical Modeling}\label{Section:mathematicalModel}

%iven the data input of section \ref{Section:Networks}, a mathematical model formulation is employed. 
Next we provide a generic, network agnostic, mixed integer nonlinear programming formulation (MINLP) for our optimisation problem. The use of each constraint set in (MINLP) depends on data availability, e.g. the demand of flow rates at a node of the network, activate Constraints (\ref{eq:nl2}). Our formulation introduces the following sets of variables: a) $x_e$,  $0\leq x_{e}\leq C_{e}$, to determine the flowrate of each edge $e\in E$, b) $c_{jp}$,  $c_{jp} \geq 0$, to determine the quality values of component $i\in I$ for pollutant $p\in P$, and c) $Y_e$,  which are binary variables equal to 1, if edge $e\in E$ is used i.e. $x_e > 0$, and 0 otherwise. In more detail, Constraints (\ref{eq:nl1}) ensure that the sum of flowrates of all edges $\delta^{j+}$ must be equal to the available provided quantity $q_{j}$, if $j$ is a flow provider and a fixed flow rate should be inserted in the system network. By (\ref{eq:nl2}), if $j$ is a flow receiver with a demand value, then the received quantity must be greater or equal with this demand value $d_{j}$. Constraints (\ref{eq:nl3}) and \ref{eq:nl4} ensure that capacity limits both on the entry and the exit of a component are respected. 

\begingroup
\tiny
\begin{flalign}
    \text{(MINLP):} & &&\notag &&\\
    &\text{optimise        }      z && \notag &&\\[1ex] 
    &\text{\underline{Flow Conservation:}} && \notag &&\\
    &\sum_{e\in \delta^{j+}}x_{e} = q_{j} && \forall j:\delta^{j-}=\text{\O}, q_j > 0 \label{eq:nl1} &&\\
    &\sum_{e\in \delta^{j-}}x_{e} \geq d_{j} && \forall j:\delta^{j+} = \text{\O}, d_{j}>0 \label{eq:nl2} &&\\
    &\sum_{e\in \delta^{j-}}x_{e} \leq C_{j} && \forall j:\delta^{j-} \neq \text{\O}, C_{j}>0 \label{eq:nl3} &&\\
    &\sum_{e\in \delta^{j+}}x_{e} \leq C_{j} && \forall j:\delta^{j+} \neq \text{\O}, C_{j}>0 \label{eq:nl4} &&\\
    &x_{e} \leq  C_{e} \cdot  Y_{e}  && \forall e \in E, C_{e}>0 \label{eq:nl5} &&\\
    & x_{e} \geq \mu \cdot Y_{e}  && \forall e \in E \label{eq:nl6} &&\\
    &\sum_{e\in \delta^{j-}}SR_{j}\cdot x_{e} = \sum_{f\in \delta^{j+}}x_{f} && \forall j:\delta^{j-}\neq \text{\O}, \delta^{j+}\neq \text{\O} \label{eq:nl7} &&\\
    &\sum_{f\in \delta^{j+}} x_{f} \leq SF_j + (1 - Y_e) \cdot M && \forall j:\delta^{j-}\neq \text{\O}, \delta^{j+}\neq \text{\O}, \forall e \in \delta^{j-} \label{eq:nl8} &&\\
    &\sum_{f\in \delta^{j+}} x_{f} \geq  SF_j - (1 - Y_e) \cdot M && \forall j:\delta^{j-}\neq \text{\O}, \delta^{j+}\neq \text{\O}, \forall e \in \delta^{j-} \label{eq:nl9} &&\\
    &\sum_{f\in \delta^{j+}} x_{f} \leq SF_j \cdot \sum_{e\in \delta^{j-}} Y_{e} && \forall j:\delta^{j-}\neq \text{\O}, \delta^{j+}\neq \text{\O} \label{eq:nl10} &&\\
    & && \notag &&\\
    &\text{\underline{Quality adjustment:}} && \notag &&\\
    &c_{jp} = \bar{c}_{jp} && \forall j: \delta^{j-} = \text{\O}, p\in P \label{eq:nl11} &&\\
    &c_{jp} = RR_{jp} \cdot \frac{\sum_{e\in \delta^{j-}}x_{e} \cdot c_{ip}}{\sum_{e\in \delta^{j-}}x_{e}} && \forall j:\delta^{j-}\neq \text{\O}, p\in P \label{eq:nl12} &&\\
    & c_{jp} \geq l_{jp} \cdot RR_{jp}  && \forall j:\delta^{j-}\neq \text{\O}, p\in P \label{eq:nl13} &&\\
    &c_{jp} \leq u_{jp} \cdot RR_{jp}  && \forall j:\delta^{j-}\neq \text{\O}, p\in P \label{eq:nl14} &&\\
    &c_{jp} - RF_{jp} \leq M \cdot (1 - Y_{e}) && \forall j:\delta^{j-}\neq \text{\O}, \forall e \in \delta^{j-}, p\in P \label{eq:nl15} &&\\
    &RF_{jp} - c_{jp} \leq M \cdot (1 - Y_{e}) && \forall j:\delta^{j-}\neq \text{\O}, \forall e \in \delta^{j-}, p\in P  \label{eq:nl16} &&\\
    %&c_{jp} - RR_{jp} \cdot c_{ip} \leq M \cdot (1 - Y_{ij}) && \forall j:\delta^{j-} = \{i\}, p\in P, \exists  RR_{jp} \label{eq:nl12} &&\\
    %&c_{ip} \cdot RR_{jp} - c_{jp} \leq M \cdot (1 - Y_{ij}) && \forall j:\delta^{j-} = \{i\}, p\in P, \exists  RR_{jp} \label{eq:nl13} &&\\
    &l_{jp} \leq \frac{\sum_{e\in \delta^{j-}}x_{e} \cdot c_{ip}}{\sum_{e\in \delta^{j-}}x_{e}} && \forall j: \delta^{j-} \neq \text{\O}, p\in P \label{eq:nl17} &&\\
    &u_{jp} \geq \frac{\sum_{e\in \delta^{j-}}x_{e} \cdot c_{ip}}{\sum_{e\in \delta^{j-}}x_{e}} && \forall j: \delta^{j-} \neq \text{\O}, p\in P \label{eq:nl18} &&\\
    & && \notag &&\\
    &\text{\underline{Variables:}} && \notag &&\\
    &0\leq x_{e}\leq C_{e} && \forall e\in E \notag &&\\
    &c_{jp} \geq 0 && \forall j\in I, p\in P \notag && \\
    & Y_e \in \{0, 1\} && \forall e\in E \notag &&
\end{flalign}
\endgroup

Constraints (\ref{eq:nl5}) and (\ref{eq:nl6}) ensure that if a flow rate traverses edge $e$, then $Y_e$ is equal to 1 (and vice versa). By (\ref{eq:nl7}), if $j$ is an intermediate node, flows are conserved, given the reduction rate $SR_{j}$. Constraints (\ref{eq:nl8}), (\ref{eq:nl9}), (\ref{eq:nl10}) refer to the case, where instead of providing the $SR_{j}$ value, the network is equipped with a component, whose outlet is always set to $SF_{j}$. More specifically, if and only if component $j$ receives any flow from at least another component, then, the sum of all the flowrates of edges leaving $j$ must be equal to $SF_{j}$.  By (\ref{eq:nl11}), if $j$ is a flow provider, then the respective quality parameter values are pre-defined. Constraints (\ref{eq:nl12}) ensure that if $j$ is a blending point (i.e., it receives flows from at least two edges), then the quality parameter values of $j$ are the weighted average of the incoming flows and the respective quality parameter values. Constraints (\ref{eq:nl13}) and (\ref{eq:nl14}), ensure that quality requirements are respected. Note, that, since $c_{jp}$ refers to the quality parameter values on the exit of component $j$, we use the $RR_{jp}$ parameter in order to calculate the respective quality value upon the component entry. On the other hand, Constraints (\ref{eq:nl15}) - (\ref{eq:nl18}) refer to the case, where the reduction on a component $j$ for a pollutant $p$ is fixed, and provided by $RF_{jp}$. As a result, Constraints (\ref{eq:nl15}) and (\ref{eq:nl16}) ensure that if component $j$ receives any amount of flowrate, then the value of $c_{jp}$ will be equal to $RF_{jp}$. Finally, Constraints (\ref{eq:nl17}) and (\ref{eq:nl18}) ensure that required levels of quality on the entry of $j$ are respected, by calculating the weighted average of concentration $c_{jp}$ and input flowrate $x_{e}$.

Note that (MINLP) is quite generic to adopt various objective functions, including fresh water minimisation, wastewater treatment maximisation, cost or energy consumption minimisation, depending on the availability of data information. For example, let $I^\prime \subseteq I$ be the set of components where cost or energy calculations are involved. Then, the following quantities correspond to the total cost and energy consumption, respectively, Cost: $\sum_{e \in \delta^{j-}} \sum_{j\in I^\prime} x_{e} \cdot VC_{j} + Y_{e} \cdot FC_{j}, VC_{j}, FC_{j} \geq 0 $, and Energy consumption: $\sum_{e \in \delta^{j-}} \sum_{j\in I^\prime} x_{e} \cdot VE_{j} + Y_{e} \cdot FE_{j}, VE_{j}, FE_{j} \geq 0$.

\subsection{Constraint Programming model}
    The MINLP includes conditional constraints, which are easily linearised using big-M constraints (e.g., (\ref{eq:nl8}) and (\ref{eq:nl9})). It also involves products of variables, necessitating more advanced linearisation techniques, as detailed in Section \ref{Section:linearApproximation}. Alternatively, the problem can be effectively addressed by constructing a Constraint Programming (CP) formulation, which efficiently handles both conditional constraints and non-linear expressions:
\begingroup
    \tiny
    \begin{flalign}
        \text{optimise}   & \,\,  z && \notag &&\\
        &\sum_{e\in \delta^{j+}}x_{e} = q_{j} && \forall j: \delta^{j-} = \text{\O}, q_{j} > 0 \label{eq:c1} &&\\
        &\sum_{e\in \delta^{j-}}x_{e}\geq d_{j} && \forall j: \delta^{j+} = \text{\O}, d_{j} > 0 \label{eq:c2} &&\\
        &\sum_{e\in \delta^{j-}}x_{e}\leq C_{j} && \forall j: \delta^{j-} \neq \text{\O}, C_{j} > 0 \label{eq:c3} &&\\
        &\sum_{e\in \delta^{j+}}x_{e}\leq C_{j} && \forall j: \delta^{j+}\neq \text{\O}, C_{j} > 0 \label{eq:c4} &&\\
        &\sum_{e\in \delta^{j-}}SR_{j}\cdot x_{e} = \sum_{f\in \delta^{j+}}x_{f} && \forall j: \delta^{j-}\neq \text{\O}, \delta^{j+}\neq \text{\O} \label{eq:c5} &&\\
        &\text{if }x_{e} > 0 \rightarrow \sum_{f\in \delta^{j+}}x_{f} = SF_{j} && \forall j: \delta^{j-}\neq \text{\O}, \delta^{j+}\neq \text{\O}, e\in \delta^{j-} \label{eq:c6} &&\\
        &c_{jp} = \bar{c}_{jp} && \forall j: \delta^{j-} = \text{\O}, p\in P \label{eq:c7} &&\\
        &c_{jp} = RR_{jp}\cdot \frac{\sum_{e\in \delta^{j-}}x_{e}\cdot c_{jp}}{\sum_{e\in \delta^{j-}}x_{e}} && \forall j: \delta^{j-}\neq \text{\O}, p
        \in P \label{eq:c8} &&\\
        &c_{jp}\leq l_{jp}\cdot RR_{jp} && \forall j: \delta^{j-}\neq \text{\O}, p\in P \label{eq:c9} &&\\
        &c_{jp}\geq u_{jp}\cdot RR_{jp} && \forall j: \delta^{j-}\neq \text{\O}, p\in P \label{eq:c10} &&\\
        &\text{if }x_{e} > 0 \rightarrow c_{jp} = RF_{jp} && \forall j\in \delta^{j-}\neq \text{\O}, e\in \delta^{j-}, p\in P \label{eq:c11} &&\\
        &l_{jp} \leq \frac{\sum_{e\in \delta^{j-}}x_{e} \cdot c_{ip}}{\sum_{e\in \delta^{j-}}x_{e}} && \forall j: \delta^{j-} \neq \text{\O}, p\in P \label{eq:c12} &&\\
        &u_{jp} \geq \frac{\sum_{e\in \delta^{j-}}x_{e} \cdot c_{ip}}{\sum_{e\in \delta^{j-}}x_{e}} && \forall j: \delta^{j-} \neq \text{\O}, p\in P \label{eq:c13} &&\\
        & && \notag &&\\
        &0\leq x_{e}\leq C_{e} && \forall e\in E \notag &&\\
        &c_{jp}\geq 0 && \forall j\in I, p\in P \notag &&
    \end{flalign}
    \endgroup
    Constraints (\ref{eq:c1}) - (\ref{eq:c4}) replace constraints (\ref{eq:nl1}) - (\ref{eq:nl4}) respectively. Constraints (\ref{eq:c5}) are equivalent with constraints (\ref{eq:nl7}). Big-M constraints (\ref{eq:nl8}) - (\ref{eq:nl10}) can be redefined, using conditional constraints, as shown in (\ref{eq:c6}). Regarding the quality adjustment constraints, it is easily seen that (\ref{eq:c7}) - (\ref{eq:c10}) are identical with (\ref{eq:nl11}) - (\ref{eq:nl14}), and (\ref{eq:c12}) - (\ref{eq:c13}) are identical with (\ref{eq:nl17}) - (\ref{eq:nl18}). Big-M constraints (\ref{eq:nl15}) and (\ref{eq:nl16}) can be reformulated as in (\ref{eq:c11}). Variables $x_{e}$ indicate the flowrates of edges $e$ and variables $c_{jp}$ indicate the quality of component $j$ for pollutant $p$.
    
    In a CP context, the variables $Y_{e}$ are unnecessary, as they can be replaced by the conditions $x_{e} > 0$, as shown in constraints (\ref{eq:c6}) and (\ref{eq:c11}). Additionally, CP allows the use of variable products. It is worth noting that variables in CP are strictly integer, so all parameters must be scaled by an appropriate power of 10, depending on the desired precision.

\subsection{Linearised model}\label{Section:linearApproximation}

\begin{figure}[h!]  
\centering
\includegraphics[scale=0.75]{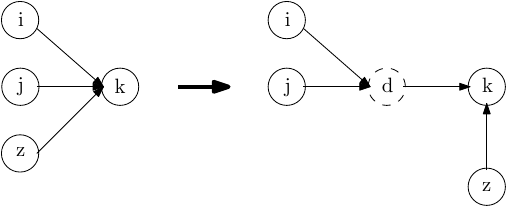}
\caption{Insertion of dummy nodes on the preprocessing stage}
\label{fig:dummy_nodes_prep}
\end{figure}

Apart from the CP formulation, it is possible to linearise the initial MINLP, by applying the linearisation proposed in \cite{Lodi15}. To linearise (MINLP), we define a discretisation number $K$ (which is typically set to 100). Consider a simple example where two nodes, $i_1$, $i_2$, are linked to a node $j$, making this a blending point with $|\delta^{j-}| = 2$. The use of discretisation allows us to calculate the parts of the mixture occupied from the edge originating from each of the nodes $i_1$ and $i_2$, i.e., let $k$ be the parts of the mixture occupied by $i_1$, then the parts occupied by $i_2$ are $K - k$. Note that, since this method is only applicable on cases where $|\delta^{j-}| = 2$, if $|\delta^{j-}| > 2$, then (on a pre-processing stage) we could create dummy intermediate nodes, and thus, transforming the blending points as if $|\delta^{j-}| = 2$, as shown on Figure \ref{fig:dummy_nodes_prep}. 
This procedure ensures that each node is the end node of no more than two other nodes; note that, by applying it iteratively, we may deal with more origin nodes on a similar manner.

Similarly with (MINLP), the linearised formulation, denoted below as (MILP), considers three sets of variables: a) $x_{e}$, $0\leq x_e\leq C_e$, that measures the flowrate of edge $e\in E$, b) $c_{jp}$, $l_{jp}\leq c_{jp}\leq u_{jp}$, which count the quality in $j\in I$ for pollutant $p\in P$, and c) the binary variables $z_{ijk}$, which are equal to 1 if the flowrate from $i$ to $j$ occupies $\frac{k}{K}$ parts of the mixture in $j$, and 0 otherwise, e.g., if $z_{ij50} = 1$ and $K = 100$, then the flow from $i$ to $j$ is equal with $50\%$ of the mixture in blending point $j$. Moreover, we introduce almost every constraint set of (MINLP), except from nonlinear Constraints (\ref{eq:nl12}), (\ref{eq:nl17}) and (\ref{eq:nl18}) and Constraints (\ref{eq:nl13}), (\ref{eq:nl14}). 

\begingroup
\tiny
\begin{flalign}
\text{(MILP):} & &&\notag &&\\
    \text{optimise}   & \,\,  z && \notag &&\\ 
    &\text{\underline{Flow Conservation:}} && \notag &&\\
    &\sum_{k = 0}^{K}z_{ijk} = 1 && \forall (i,j) = e^{1}, \delta^{j-} = \{e^{1}, e^{2}\} \label{eq:l1} &&\\
    &z_{ijk} = z_{rj(K-k)} && \forall \delta^{j-} = \{(i,j), (r,j)\}, k = 1, ..., K - 1 \label{eq:l22} &&\\
    &(K - k)\cdot x_{e^1} - k\cdot x_{e^2} \leq M\cdot (1 - z_{ijk}) && \forall (i,j) = e^{1}, \delta^{j-} = \{e^{1}, e^{2}\}, k = 1, ..., K - 1 \label{eq:l2} &&\\
    &k\cdot x_{e^2} - (K - k)\cdot x_{e^1} \leq M\cdot (1 - z_{ijk}) && \forall (i,j) = e^{1}, \delta^{j-} = \{e^{1}, e^{2}\}, k = 1, ..., K - 1 \label{eq:l3} &&\\
    &z_{ij0} + Y_{ij}  = 1 && \forall \delta^{j-} = \{(i,j), (r,j)\} \label{eq:l4} &&\\
    &z_{ijK} + Y_{rj}  \leq 1 && \forall \delta^{j-} = \{(i,j), (r,j)\} \label{eq:l5} &&\\
    & && \notag &&\\
    &\text{\underline{Quality adjustment:}} && \notag &&\\
    &c_{ip} \cdot RR_{jp} -   c_{jp} \leq M \cdot (1 - Y_{ij}) && \forall j:\delta^{j-} = \{(i,j)\}, |\delta^{j-}| = 1,  p\in P \label{eq:l6} &&\\
    &c_{jp} - c_{ip} \cdot RR_{jp}\leq M \cdot (1 - Y_{ij}) && \forall j:\delta^{j-} = \{(i,j)\}, |\delta^{j-}| = 1, p\in P \label{eq:l7} &&\\
    & RR_{jp} \cdot \frac{k\cdot c_{ip} + (K - k)\cdot c_{rp}}{K} - c_{jp}  \leq M\cdot (1 - z_{ijk}) && \forall \delta^{j-} = \{(i,j), (r,j)\}, k = 1, ..., K - 1, p\in P \label{eq:l8} &&\\
    &c_{jp} - RR_{jp} \cdot \frac{k\cdot c_{ip} + (K - k)\cdot c_{rp}}{K}\leq M\cdot (1 - z_{ijk}) && \forall \delta^{j-} = \{(i,j), (r,j)\}, k = 1, ..., K - 1, p\in P \label{eq:l9} &&\\
    &c_{rp} \cdot RR_{jp} - c_{jp} \leq M\cdot (2 - z_{ij0} - Y_{rj}) && \forall \delta^{j-} = \{(i,j), (r,j)\}, p\in P \label{eq:l10} &&\\
    &c_{jp} - c_{rp} \cdot RR_{jp} \leq M\cdot (2 - z_{ij0} - Y_{rj}) && \forall \delta^{j-} = \{(i,j), (r,j)\}, p\in P \label{eq:l11} &&\\
    &c_{ip} \cdot RR_{jp} - c_{jp} \leq M\cdot (1 - z_{ijK}) && \forall \delta^{j-} = \{(i,j), (r,j)\}, p\in P \label{eq:l12} &&\\
    &c_{jp} - c_{ip} \cdot RR_{jp} \leq M\cdot (1 - z_{ijK}) && \forall \delta^{j-} = \{(i,j), (r,j)\}, p\in P \label{eq:l13} &&\\
    &c_{jp} \geq l_{jp} && \forall \delta^{j-} = \{(i,j)\}, p\in P \label{eq:l14} &&\\
    &c_{jp} \leq u_{jp} && \forall \delta^{j-} = \{(i,j)\}, p\in P \label{eq:l15} &&\\
    &l_{jp} \leq \frac{k\cdot c_{ip} + (K - k)\cdot c_{rp}}{K} +  M\cdot (1 - z_{ijk}) && \forall \delta^{j-} = \{(i,j), (r,j)\}, k = 1, ..., K - 1, p\in P \label{eq:l16} &&\\
    &u_{jp} \geq \frac{k\cdot c_{ip} + (K - k)\cdot c_{rp}}{K} -  M\cdot (1 - z_{ijk}) && \forall \delta^{j-} = \{(i,j), (r,j)\}, k = 1, ..., K - 1, p\in P \label{eq:l17} &&\\
    &l_{jp} \leq c_{rp} + M\cdot (2 - z_{ij0} - Y_{rj}) && \forall \delta^{j-} = \{(i,j), (r,j)\}, p\in P \label{eq:l18} &&\\
    &l_{jp} \leq c_{ip}  + M\cdot (1 - z_{ijK}) && \forall \delta^{j-} = \{(i,j), (r,j)\}, p\in P \label{eq:l19} &&\\
    &u_{jp} \geq c_{rp} - M\cdot (2 - z_{ij0} - Y_{rj}) && \forall \delta^{j-} = \{(i,j), (r,j)\}, p\in P \label{eq:l20} &&\\
    &u_{jp} \geq c_{ip} - M\cdot (1 - z_{ijK}) && \forall \delta^{j-} = \{(i,j), (r,j)\}, p\in P \label{eq:l21} &&\\
    & && \notag &&\\
    &\text{\underline{Variables:}} && \notag &&\\
    &0\leq x_{e}\leq C_{e} && \forall e\in E \notag &&\\
    &c_{jp} \geq 0 && \forall j\in I, p\in P \notag &&\\
    &Y_e\in \{0, 1\} && \forall e\in E \notag &&\\
    &z_{ijk}\in \{0, 1\} && \forall (i,j) = e^{1}, \delta^{j-} = \{e^{1}, e^{2}\}, k = 0, ..., K \notag &&
\end{flalign}
\endgroup

In more detail, Constraints (\ref{eq:l1}) on each blending point ($|\delta^{j-}| = 2$) are mandatory to assign at least one $k$ value to the respective edge $(i,j)$ (that value could be also 0 or $K$). Constraints (\ref{eq:l22}), ensure that if any flow rate is drawn from both $i$ and $r$, on a blending point $j$, and if the compartment from $i$ is equal to $k$, then the respective compartment from $r$ is equal to $K-k$.  Furthermore, if $j$ is a blending point (i.e., it receives flows from two edges $e^{1}$, $e^{2}$), then edge $e^{1}$ covers $k$ parts out of $K$ in the mixture (\ref{eq:l2}). If $e^{1}$ covers $k$ parts, then $e^{2}$ should cover the remaining $K - k$, by Constraints (\ref{eq:l3}). Constraints (\ref{eq:l4}) ensure that, if $i$ occupies 0 counterparts on the mixture (i.e., $z_{ij0} = 1$), then no flowrate is traversing edge $ij$ and thus, $Y_{ij} = 0$. Similarly, Constraints (\ref{eq:l5}) ensure that if flowrate from component $i$ occupies all the mixture ($z_{ijK}$ = 1), then no flowrate is drawn from component $r$ ($Y_{rj} = 0$). Furthermore, it ensures that it is not mandatory to draw flowrates from at least one of the components $i$, $r$, as both variables may be equal to 0. Constraints (\ref{eq:l6}) - (\ref{eq:l13}) refer to the combination of components $j$ and pollutants $p$, which are associated with the $RR_{jp}$ attribute. By (\ref{eq:l6}) and (\ref{eq:l7}), it is ensured that if a component $j$ receives flowrates from exactly one component $i$, then $c_{jp}$ will be equal to the product of $RR_{jp}$ and $c_{ip}$. In addition, if $e^{1}$ covers $k$ parts, then $c_{jp}$ is determined by the mean average formula, by constraints (\ref{eq:l8}), (\ref{eq:l9}). However, this is only applied when flowrates are drawn from both $e^{1}$ and $e^{2}$. If flowrates are drawn only from $r$, then by Constraints (\ref{eq:l10}),  (\ref{eq:l11}) it is ensured that $c_{jp}$ is equal to $c_{rp}$. Similarly, Constraints (\ref{eq:l12}), (\ref{eq:l13}) ensure the opposite scenario (flowrates only drawn from $i$). Constraints (\ref{eq:l14}) - (\ref{eq:l21}) correspond to the case where $RF_{jp}$ is provided. Constraints (\ref{eq:l14}) and (\ref{eq:l15}) ensure that the quality limits are respected when such a component receives flowrates from exactly on component. Similarly, Constraints (\ref{eq:l16}) - (\ref{eq:l21}) enforce the quality requirements when component $j$ draws flowrates from two components, $i$ and $r$.

\section{DSS analysis, design and implementation} \label{Section:DSSAnalysis}

Having presented the detailed technical aspects of the mathematical formulation, in this section, we provide a detailed analysis of the proposed DSS. To do this, the Unified Modelling Language (UML) is employed to analyze the user requirements (Section \ref{Section:user_requirements}), and the respective data requirements (Section \ref{Section:data_model}). Last, we demonstrate through a component diagram the architecture of the DSS (Section \ref{Section:system_architecture}) that realises the user requirements.  

\subsection{User requirements} \label{Section:user_requirements}
Based on the case studies in Section \ref{Section:CaseStudies}, we identified key functionalities for end-users to design or operate water networks in industrial settings. For network design, users configure the network structure with essential components (e.g., pipes, pumps, tanks, cleaning systems, valves) and set operational parameters (e.g., flow limits, water quality metrics). If the network is still under investigation rather than implemented, users need tools to assess design decisions and evaluate its performance. For network operation, the focus shifts to optimizing the functionality of an existing network, allowing users to configure its structure and define objectives such as maximizing water reuse.

\begin{figure}[h!]
\centering
\includegraphics[scale=.35]{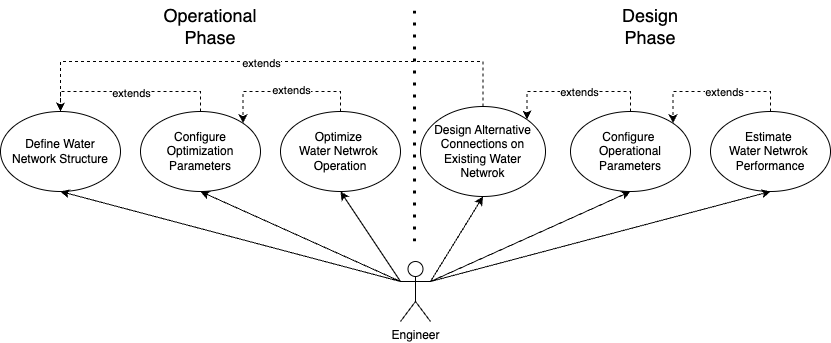}
\caption{Use Case Diagram}
\label{fig:use-case-diagram}
\end{figure}

Hence, the functional requirements of the proposed DSS involve six use cases (UCs) grouped in two categories: (a) UCs related to the operational phase of the water network and (b) UCs related to the design phase of the water network. Our goal is to create a single system that accommodates both types of functional requirements, making it essential to distinguish between the two cases. This categorization benefits system developers by clarifying different system use cases, while also helping end-users clearly differentiate between design and operational functionalities. However, this approach may be less useful for organizations that require only one of the two functionalities, rendering the other category redundant.
Figure \ref{fig:use-case-diagram} shows these UCs with the categories and the actor, i.e., the Engineer. We opt for a generic enough actor type to accommodate all the involved organisations, which usually employ chemical or production engineers to deal with such decisions.

Regarding the Operational Phase category, the ``Define Water Network" UC relates to the actions required to set up the actual structure of the water network. This involves the depiction of the network as a graph, where the nodes are a network component (e.g., a tank, a water-treatment unit, a production process, etc) and the edges of graph indicate the connections (i.e., water pipes) between two components. An edge may act as an incoming/outgoing connection to/from a component, while a node may have more than one edge as input or output. In addition, every node may include specs on the water use. For instance, a water treatment unit may intake water of a maximum volume  with specific contaminant levels on its input and output. At any time, the user wishes to edit these parameters and update the water network configuration. 

The ``Configure Optimisation Parameters" UC relates to the actions required to set up the required optimisation parameters, objectives and KPIs. In more detail, the user could set different objectives considering that there are different decisions to be accommodated. For instance, this could be the maximisation of the treated water or minimisation of the fresh-water intake on a specific set of edges. In addition, the user may set specific parameters (e.g., contaminants levels) on a set of nodes. Last, the "Optimize Water Network Operation" UC relates to the actions required to trigger the optimisation algorithm and review the results. 

Regarding the Design Phase category, the ``Design Alternative Connections on Existing Network Structure" UC relates to the actions required to define new network components or non-existing connections between components in order to asses the performance of a new network structure. As previously, the definition of a new connection involves the generation of a new edge between two components, including the related parameters that affect the decision (e.g., capacity, contaminants, etc). The ``Configure Operational Parameters" UC relates to actions required to configure the parameters of the water network before assessing its performance. This includes the definition of the parameters domain (i.e., single value or range), the distribution that rules a parameter defined within a range or the iterations that the optimisation method will run. Last, the ``Estimate Water Network Performance" relates to the actions required to evaluate the performance of the water network based on the configuration of operational parameters in the previous UC. If the user has defined ranges for a set of parameters with the respective distributions, then it is required to include within the results the generated instances that were used for performance estimation input, while also indicate the most dominant instances.

Concerning non-functional requirements identified in the case studies, the system must effectively handle a considerable amount of data. Security concerns are paramount. Access to the system is contingent upon successful authentication, and exporting data without authorisation is strictly prohibited. Additionally, the system must be compatible with various hardware and operating systems. Both vertical and horizontal scalability are essential. Access to the system must be available 24/7 with near-perfect reliability. Last, all capabilities of the system and its results should be available and easy to integrate with existing systems.

\subsection{Data requirements} \label{Section:data_model}
\begin{figure}[h!]    
\begin{center}
\includegraphics[scale=0.37]{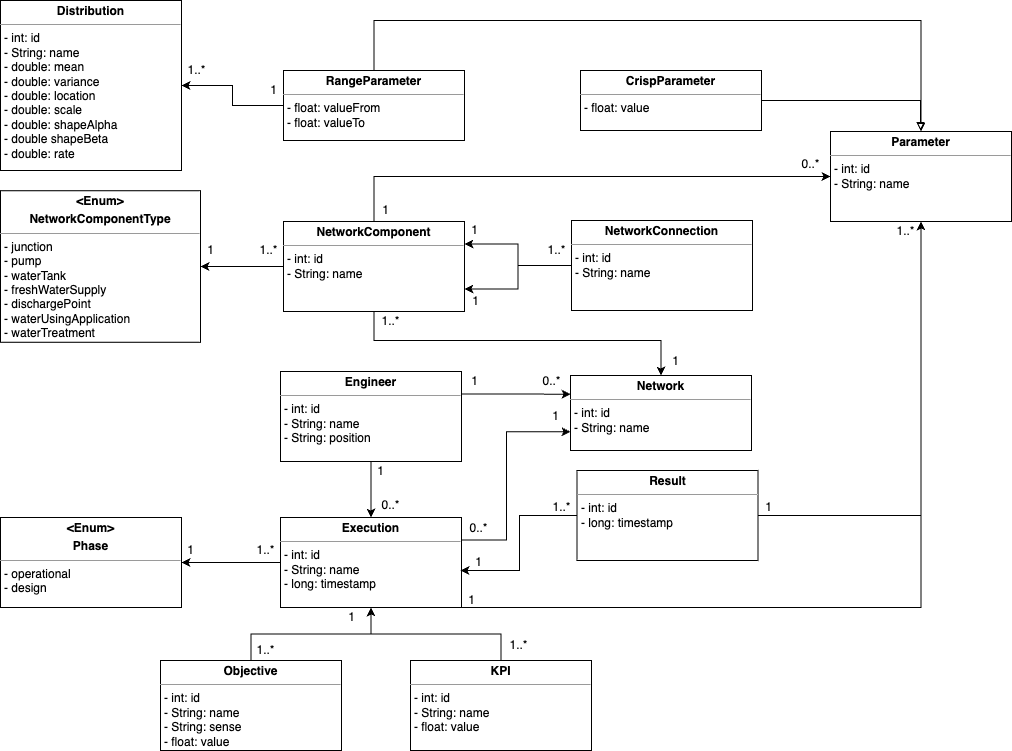}
\end{center}
\caption{Domain Diagram}
\label{fig:domain-diagram}
\end{figure}

In this section, we describe the data requirements in the form of a domain diagram (Figure \ref{fig:domain-diagram}). This way, we provide a high-level overview of the domain under consideration and the related data requirements that can complement and enable the implementation of the user requirements.

The central entity of the domain model is the \textbf{Network}. It depicts the water network of the industrial site following the structure of a graph. A Network is composed of one or more \textbf{NetworkComponents} which refer to the nodes of the graph. These are of a specific \textbf{NetworkComponentType} which can be a water treatment station, a water tank, fresh water supply, etc. Every node, may have specific \textbf{Parameters} referring to the water specifications of its input and output. These parameters can be either \textbf{CrispParameters}, i.e., clear-cut values that are not subject to interpretation or uncertainty or \textbf{RangeParameters}, i.e., parameters that are defined within a specific continuous range of values. In the latter case, the parameters can be associated with a specific \textbf{Distribution}, i.e., the description of how the probability is spread out over the given parameter range.  A \textbf{NetworkComponent} can be connected with another via a \textbf{NetworkConnection}, i.e., the edges of the water network graph. This concludes the description of the \textbf{Network} as a whole, which is usually modeled by the end-user of the system, i.e., the \textbf{Engineer}. 

On the decision support part of this model, the \textbf{Execution} entity depicts the possible runs of the involved optimisation method. Such executions refer to the operational phase, or the design phase, while both of them are associated with one ore more \textbf{Objectives} and a set of \textbf{KPIs}. In addition, they can be associated with specific \textbf{Parameters} as input, while they provide \textbf{Results} which again contain a set of parameters with their respective values as identified for each case accordingly. 

\subsection{System Architecture} \label{Section:system_architecture}
\begin{figure}[h!]    
\centering
\includegraphics[scale=0.39]{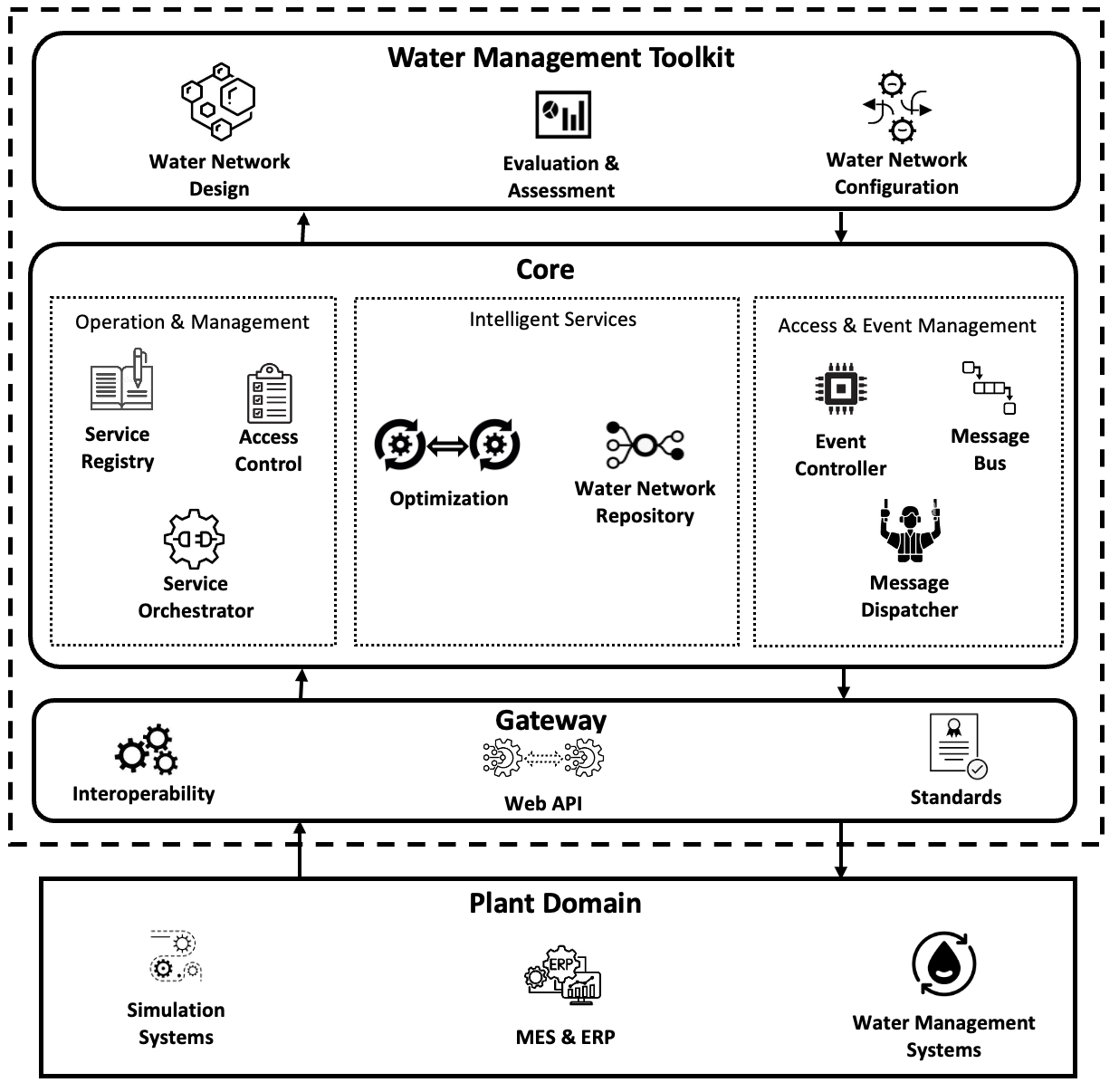}
\caption{DSS architecture}
\label{fig:dss-architecture}
\end{figure}

Figure \ref{fig:dss-architecture} shows the architecture of the proposed DSS. The system is comprised of three layers (as seen in the hyphenated area): the \textbf{Water Management Toolkit}, the \textbf{Core} and the \textbf{Gateway}. The \textbf{Plant Domain} area refers to external systems that exist in the premises of a plant, such as a Manufacturing Execution System (MES) or Enterprise Resource Planing system (ERP), a Water Management System or Simulation Systems that are being used for production purposes. These systems may act both as data sources or output receivers for/from the proposed DSS. Let us now describe the details of the DSS components. 

Starting from the bottom of the DSS, the \textbf{Gateway} layer is the entry point and output source of data flowing in or out of the DSS respectively. This layer works as an intermediate between the existing systems in the plant premises and the core of the DSS. It is comprised of three components. The \textbf{Interoperability} component includes the implementation of interfaces for integrating external systems. The \textbf{Web API} component is the main entry point of incoming data and requests from external systems, as well as an output point of results from the DSS. The \textbf{Standards} component is a repository for interfaces and standards utilised for all communication either internally or externally.

In the middle, lies the \textbf{Core} of the DSS, which is comprised of three sub-elements. The first, i.e., \textbf{Operations \& Management} sub-element contains three components. The \textbf{Service Registry} component is a central record of all the employed services, along with a resource catalogue available to any authorised user or service. The \textbf{Access Control} component is responsible for granting access to services or users of the system. The \textbf{Service Orchestrator} component manages the internal operations of the system. Any service required to be (re-) deployed is handled by this module. It also handles the (re)allocation of available system resources. Second, the \textbf{Intelligent Services} sub-element, is comprised of two components. The \textbf{Optimisation} component is responsible for the optimisation of the water network use in its operational and design phase. The \textbf{Water Network Repository} component is responsible for keeping the current water network structure as well as designed/reconfigured ones (which are under consideration), which enables the assessment of different design or operational scenarios. Third, the \textbf{Intelligent Services} sub-element, is comprised of three components. The \textbf{Event Controller}  component is responsible to process incoming events and trigger the respective system component(s). The \textbf{Message Bus} component is the main and central communication mechanism for all system components which transfers data and information in a scalable and asynchronous manner. The \textbf{Message Dispatcher} component is responsible for dispatching messages from one component to another standing on top of the message bus.

At the top of the DSS lies the \textbf{Water Management Toolkit} element, which refers to the user interface that basically covers three aspects of functionality via three components: i) \textit{Water Network Design}, that enables the design of a water network structure, ii) \textit{Evaluation \& Assessment}, that provides useful dashboards to visualise optimisation results, and iii) \textit{Water Network Configuration}, that enables the configuration of the parameters of a water network structure.

\section{DSS usage scenarios}

\label{Section:DSSUseScenaria}
In this section, the methodology designed is applied on three case studies. We begin by providing the experimental setting under which the proposed DSS was tested (Section \ref{Section:experimentalMethodology}). Then, we describe in detail the usage scenario for each case study as described in Section \ref{Section:CaseStudies}, along with the results of the DSS for each case study (Sections \ref{Section:impactOilRefinery}, \ref{Section:impactChemicalA} and \ref{Section:impactChemicalB}).

\subsection{Experimental Setting} \label{Section:experimentalMethodology}
Experiments are performed on three phases: \textbf{design}, \textbf{operational} and \textbf{comparison}. On the design phase experiments, all possible decisions are incorporated as input data to the DSS. Then, fields \textbf{optionsCompared} are employed, where the candidate options are included. Due to different possible parameter values, e.g., concentration of a quality parameter on a stream, we apply a flexible approach by supplying the algorithm with parameter ranges (lower and upper limits) and distributions (e.g., uniform) in order to determine the most appropriate candidate option under different input scenarios. For a fixed number of trials $N$, an input dataset is generated based on the parameters' ranges and distributions, and then, the optimisation model is executed. When the solution is returned, the flowrates of the used edges (i.e. the ones with flow rate $> 0$) are adjusted thus, obtaining the candidate option, which was found in the corresponding solution. The above procedure is shown in Algorithm \ref{alg:alg_design}. 

\begin{algorithm}[h]
	\fontsize{9}{11}\selectfont
	\caption{ \label{alg:alg_design} \textsc{Design Phase Experiments Procedure} }
	\begin{algorithmic}[1]
	\STATE Create input set $A$ with ranges and distributions;
        \STATE Let $N$ be the number of Trials;
        \STATE Let $Res$ be a counter for each candidate option;
	\FOR {$n \in N$}
	\STATE Generate dataset $B$ based on $A$;
	\STATE Solve the optimisation model $Opt$ for input $B$;
        \STATE Let $k$ be the solution returned by $Opt$;
        \STATE Find the option used in $k$;
	\STATE Update $Res$;
  	\ENDFOR\\
	\RETURN $Res$;
	\end{algorithmic}

\end{algorithm}

Regarding the \textbf{operational phase experiments}, we employ the options selected from the design phase and aim at estimating the computational (\textit{Gap}) and time efficiency (\textit{Time}) of our approach. As a result, we, also, experiment with different values of the discretisation parameter $K$ to estimate the computational burden imposed as a trade off with greater accuracy. The cut-off for \textit{Gap} is 1 \% and for \textit{Time} 90 seconds.

In the final set of \textbf{comparison} experiments, we aim at comparing predefined metrics and KPIs of the network before and after the implementation of the decisions derived by the use of our DSS. The KPIs employed are different on each case study and have been identified as critical after discussions with the corresponding stakeholders.

All value ranges of quality parameters $(c_{ip})$, quality requirements ($l_{ip}$, $u_{ip}$), fixed and variable reduction flowrates ($SF_{i}$, $SR_{i}$, $RF_{ip}$ and $RR_{ip}$) are generated based on information provided by the case studies. For convenience and due to lack of information concerning the possible parameter distribution, we assume that parameters are selected uniformly at random from their corresponding value range; yet, other distributions could be supported as well, if needed. The value of $K$ is set to 200 (unless otherwise stated) and the number of trials $N$ performed for each case study is 500. Interestingly, although $N$ is quite large, our approach appears to be highly efficient, obtaining optimal solutions within a few seconds. In terms of the number of experiments, during the \textbf{design phase}, we conduct 500 experiments per case study, resulting in a total of 1500 executions. In the \textbf{operational phase}, we focus on the Oil Refinery and Chemical Industry A cases. As outlined in Sections \ref{Section:impactOilRefinery} and \ref{Section:impactChemicalA}, we generate 500 independent instances for each case study. Since experimentation is carried out under four different values of the discretisation parameter $K$, the total number of executions per case study amounts to 2000. Additionally, the \textbf{comparison phase}, which involves solving the \textit{current network}, requires another 500 experiments per case study. Therefore, the total number of executions is calculated as: $(500 \times 3) + (500 \times 4) \times 2 + (500 \times 2) = 6,500$.

All experiments have been performed on a server with 4 Intel(R) Xeon(R) E-2126G @ 3.30GHz processors and 11 GB RAM, running CentOS/Linux 7.0. Python 3.12.0 was used for scripting and Pyomo 6.6.2. The solver of the MILP model was Gurobi Optimizer 9.1.5 \cite{gurobi}, via the Pyomo library \cite{bynum2021pyomo, hart2011pyomo}, compatible for Python 3.12.0. 

Before delving into the specific experiments for each case study, we present a key result from our initial experimentation. Using the same datasets, we found that even with a 10-minute time limit (noting that the MILP is executed in under a minute), the CP model performed significantly worse than the linearised MILP. Out of 500 instances, the CP method failed to produce a feasible solution within the time limit for 144 cases. For the remaining 356 instances, the optimality gap ranged from 62.8\% to 78.32\%, with an average gap of approximately 71.38\%. Given these results, we have decided not to present the full set of CP model outcomes for each case study.

\subsection{Oil refinery} \label{Section:impactOilRefinery}
As stated on Subsection \ref{Section:CS_Descriptions}, the network of this oil refinery contains two parallel wastewater tanks, namely T1 and T2, of different capacities and quality parameter concentrations. Flowrates drawn from each Tank are mixed before entering Treatment process Tr\_1 and finally, discharged with respect to specific limits. In case these limits are not met, the plant is facing two risks: i) severe legislation sanctions, or ii) the need to employ additional high cost chemicals in order to improve water quality. The optimisation objective is the maximisation of wastewater treatment on the whole network. Since our case study aims at dealing permanently with these risks, the designed DSS is deployed on all phases.

\begin{figure}[!htb]
   \begin{minipage}[b]{0.48\textwidth}
     \centering
     \includegraphics[scale=.44]{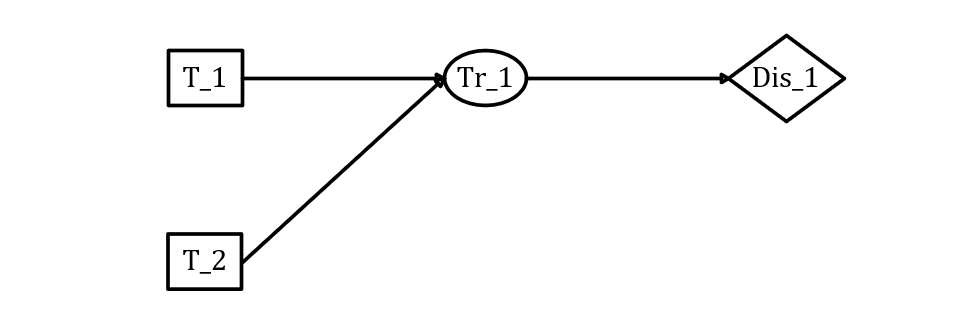}
     \caption{Current network status of Oil Refinery case study}\label{fig:case_1_as_is}
   \end{minipage}\hfill
   \begin{minipage}[b]{0.48\textwidth}
     \centering
     \includegraphics[scale=.42]{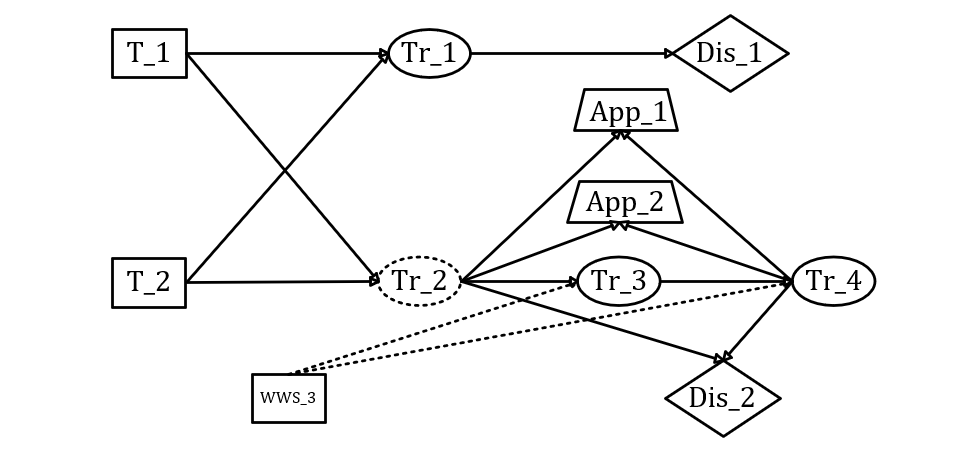}
     \caption{Updated network status of Oil Refinery case study}\label{fig:case_1_to_be}
   \end{minipage}
\end{figure}

%\begin{figure}[h!]    
%\centering
%\includegraphics[width=0.7\textwidth]{case_1_design_as_is.png}
%\caption{Current network status of Oil Refinery case study}
%\label{fig:case_1_as_is}
%\end{figure}

Under the design phase, and after consultation with inbound engineers, the insertion of more components has been recommended. The `updated' network includes: i) one wastewater stream, WWS\_3 of different composition to the tanks, ii) two applications, App\_1 and App\_2, of different requirements, where treated water should be used, iii) one already existing treatment process, Tr\_3 of known treating capability and requirements due to sensitive in-process components and finally, iv) a variable treatment process, Tr\_2, which corresponds to an investment decision with three candidate settings. \textit{Setting 1} is the variation of Tr\_2 with the greater treatment efficiency, but lower in terms of capacity and water losses, while \textit{Setting 2} provides water of worse quality, however preserving higher flowrates due to the lower degree of treatment. Finally, \textit{Setting 3} provides the most enhanced water quality  compared to the other Settings - yet it is accompanied with lower capacity and quantity preserved. The second design decision is the connection of WWS\_3 to the most appropriate component. On this scenario, there are two different options; Option A, which connects WWS\_3 to Tr\_3, aiming at the treatment of this stream, and Option B, that connects WWS\_3 to Tr\_4 (representing a Tank), where WWS\_3 will be mixed with treated water from other streams.

%\begin{figure}[h!]    
%\centering
%\includegraphics[width=0.7\textwidth]{case_1_design_to_be.png}
%\caption{Updated network status of Oil Refinery case study}
%\label{fig:case_1_to_be}
%\end{figure}

%Under the operational phase, we will conduct experiments on different possible values of input parameters (i.e. quality parameters and water flows) on the already established network to validate the capability of our system. 

%These different options are included on both the \textbf{optionsCompared} field and also, on the \textbf{conflicts} field. More specifically, on the latter we include two different constraints: a) only one of the candidate settings for Tr\_2 and b) only one of the two candidate connections of WWS\_3 will be used. As mentioned above, 

The number of trials $N$ is set to 500 and data are generated as described in Algorithm \ref{alg:alg_design}. Table \ref{tab:oil-res-design} presents the frequency of each option on the optimal solutions returned. The average execution of each trial is approximately 11 seconds, including the preprocessing stages.

\begin{table}[h!]
\centering
 \begin{tabular}{||c || c c c||} 
 \hline
  & Setting 1 & Setting 2&  Setting 3 \\ [0.5ex] 
 \hline\hline
 Option A & 209 & 7  & 3 \\ 
 Option B & 240 &  27 & 14 \\
 \hline
 \end{tabular}
 \caption{Design phase results for Oil Refinery}
  \label{tab:oil-res-design}
\end{table}

From the design phase results (Table \ref{tab:oil-res-design}), it is obvious that Option B is relatively more preferable than Option A. From an engineering point of view, it seems that it is not always possible to treat WWS\_3 in Tr\_3 due to the quality requirements of the latter. On the other hand, when mixed with already treated water, the flowrate of WWS\_3 may be diverted to applications App\_1 and App\_2 or be discharged within the appropriate quality limits. Setting 1 seems to be the dominant one with presence on the solution on the 449 out of 500 different trials. As a result, the most suitable suggestion seems to be the combination of Setting 1 with Option B, even though the differences with Option A are not dramatically large. 

Based on the above design decisions, we move on to the operational phase results, where for different values of $K$, the scalability of our method is tested, in terms of computational time, number of constraints and variables. More precisely, we consider a set of four representative values for $K\in \{100, 200, 500, 1000\}$. In order to avoid computational implications caused by different input parameters, we generate a single dataset of 500 independent input instances and test the different values of $K$ on datasets with the same parameters. Computational results are shown in Table \ref{tab:oil-res-oper}. Columns Vars, BVars and CVars correspond to the number of variables per type (B,C refer to Binary, Continuous, respectively) generated for each different value of $K$. Constraints refers to the number of constraints needed, Time (s) is the computational time in seconds and Gap refers to the optimality gap. 

\begin{table}[h!]
\centering
 \begin{tabular}{||c || c c c c c c||} 
 \hline
  $K$ & Vars & BVars &  CVars & Constraints & Time (s) & Gap (\%)\\ [0.5ex] 
 \hline\hline
 100 & 1306 & 1226 & 80 & 9211 & 2.81 & $<$ 0.3 \\ 
 200 & 2506 & 2426 & 80 & 18211 & 4.35 &  $<$ 0.3 \\
  500 & 6106 & 6026 & 80 &  45211 & 12.96   & $<$ 0.3 \\ 
 1000 & 12106 & 12026 & 80 & 90211 & 29.61 & $<$ 0.2\\
 \hline
 \end{tabular}
 \caption{Operational results on Oil Refinery}
  \label{tab:oil-res-oper}
\end{table}

CVars remain 80 at all instances, as the number of $K$ does not affect the continuous variables $c_{ip}$ and $x_e$, which only depend to the number of edges, components and quality parameters. On the other hand, there is significant increase on the number of binary variables and constraints. Interestingly, our mathematical model is highly efficient, as all instances are solved to optimality in a few seconds. However, although the solutions on every instance (for each value of $K$) are different, their objective values do not vary significantly. Thus, the examination of larger values of $K$ is interesting only in terms of greater flowrate accuracy.

The final step for our approach on the Oil Refinery case study is the comparison experiments between the current (Figure \ref{fig:case_1_as_is}) and the updated network status (Figure \ref{fig:case_1_to_be}). More precisely, we aim at examining the improvements occurred by the implementation of the updated network. As a result and similarly to the operational phase experiments, a dataset composed of 500 different input parameters is created and then executed twice: i) for the current network and ii) for the updated network. Then, the obtained solutions are evaluated in terms of the following representative KPIs: i) {\it Feasibility}, whether any amount of wastewater has been treated, ii) {\it Wastewater treated}, the average amount of wastewater which has been treated in the network, iii) {\it Amount Discharged}, the amount of wastewater which has been discharged within the specified limits, iv) {\it Amount Reused}, the amount of wastewater which has been reused, v) {\it Water Losses}, the amount of wastewater lost on treatment processes. 

\begin{table}[h!]
\centering
 \begin{tabular}{||c || c c c c c ||} 
 \hline
  $Network$ & Feas. & WW Tr. & Disch. \% & Reused \% & Losses \% \\ [0.5ex] 
 \hline\hline
 Current & 242 & 121.7 & 89.86 &  0 & 10.14 \\ 
 Updated & 470 & 991.35 & 7.05  & 89.78 & 3.17 \\
 \hline
 \end{tabular}
 \caption{Comparison of results on Oil Refinery}
  \label{tab:oil-res-comp}
\end{table}

Obtained results (Table \ref{tab:oil-res-comp}) reveal that our approach has significant impact in terms of wastewater treatment and water reuse potential. In more detail, it is noteworthy that on the current network any amount of wastewater is treated on only 242 instances (48.4 \%), while the corresponding number on the updated network stands at 470 (94 \%). In terms of the examined KPIs, a significant increase on the average flowrate of treated wastewater is observed (121.7 against 991.35). Moreover, the percentage of the discharged and the lost water significantly drops from 89.86 \% to 7.05 \% and from 10.11 \% to 3.17 \%, respectively. Finally, the most important result is highlighted on the \textit{Reused \%} column, where it is shown that moving from zero water reuse, we achieve a percentage of 89.78 \% on the new network. 

From the above findings, we may conclude that the employment of the DSS on the oil refinery case study leads to successful network design and operational decisions. The proposed DSS was able to provide numerical consultation to the plant stakeholders in order to decide upon the most efficient design possibilities. Regarding the operational phase, our experiments underline the capability of our mathematical model both in terms of solution quality and computational time efficiency. Finally, the comparison experiments on the states before and after the DSS implementation demonstrate immense improvements on the critical KPIs, as designated from the plant stakeholders.

\subsection{Chemical industry A} \label{Section:impactChemicalA}

\begin{figure}[!htb]
   \begin{minipage}[b]{0.47\textwidth}
     \centering
     \includegraphics[scale=.33]{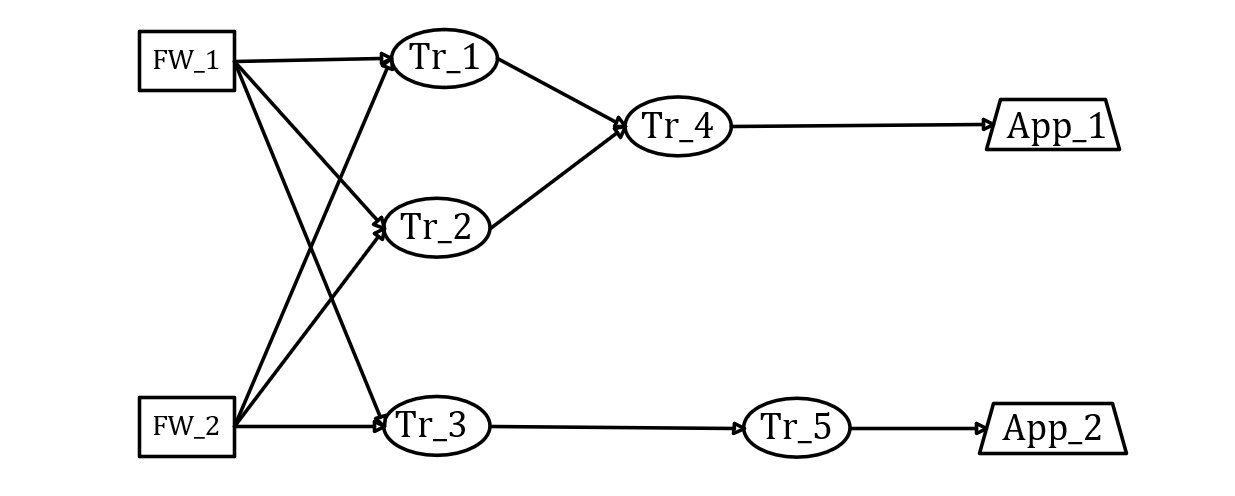}
     \caption{Current network status of Chemical Industry A }\label{fig:case_3_as_is}
   \end{minipage}\hfill
   \begin{minipage}[b]{0.45\textwidth}
     \centering
     \includegraphics[scale=.33]{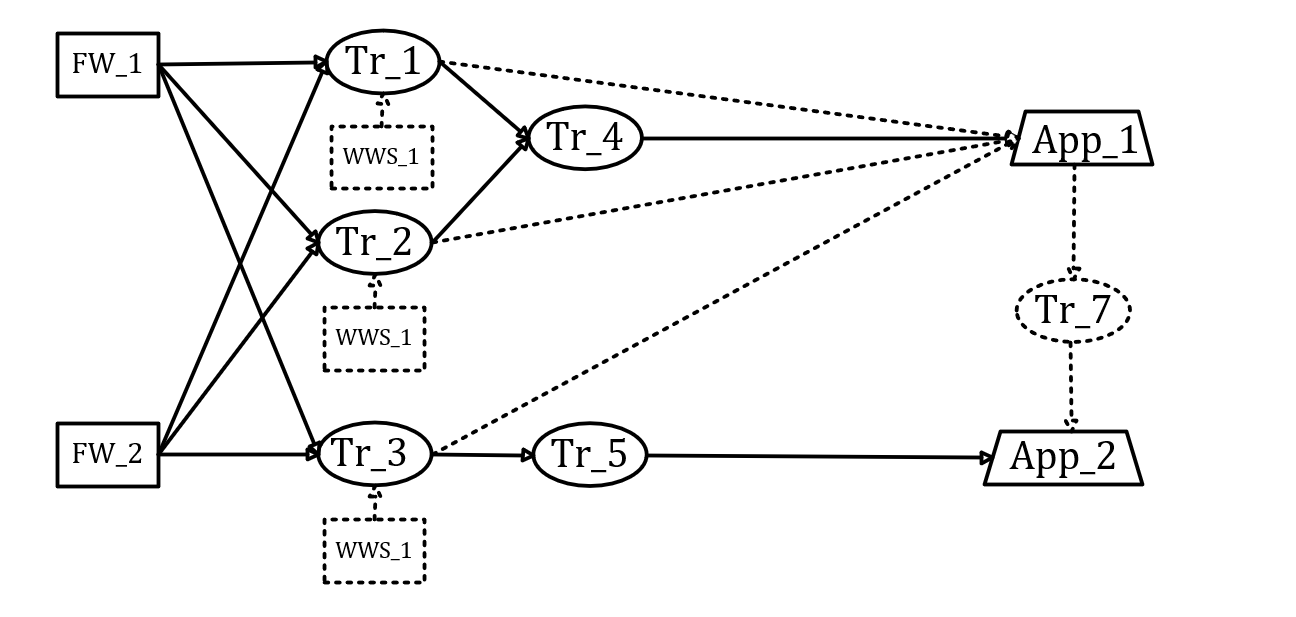}
     \caption{Updated network status of Chemical Industry A }\label{fig:case_3_to_be}
   \end{minipage}
\end{figure}

 %let us provide a description in terms of the current situation and then, present out approach. 
The initial network of Chemical Industry A is presented on Figure \ref{fig:case_3_as_is}. There are two fresh water intake sources (namely $FW_1$ and $FW_2$), which correspond to a nearby lake and a river. These two sources contribute to the provision of fresh water used for production processes and correspond to different values of quality parameters. The network consists of 2 types of treatment processes (Tr\_1, Tr\_2, Tr\_3 and Tr\_4, Tr\_5) of different  treatment capabilities and capacities. Tr\_1, Tr\_2 and Tr\_3 are connected to the fresh water sources, with the capability of mixing flowrates before the treatment. Applications, App\_1 and App\_2, are associated with different quality requirements and water demand. The optimisation goal of this case study is to minimise the freshwater intake. The described network is presented on Figure \ref{fig:case_3_as_is}.

The design goal comes from suggestions of the inbound engineers. A wastewater stream, namely WWS\_1, needs to be inserted in the system. According to expert opinion, the candidate nodes correspond to treatment processes Tr\_1, Tr\_2 and Tr\_3, as they are more suitable for the treatment of WWS\_1. The first decision lies on the exact treatment process, to which this wastewater stream will be connected. The second design decision is related to the connection of App\_1 with more components on the network. More precisely, it is examined whether App\_1 may be fed with any flowrate from Tr\_1, Tr\_2 or Tr\_3, skipping Tr\_4. This avoids some of the water losses on Tr\_4 and thus, reduce the fresh water intake demanded. The final decision targets the effluent of App\_1, which (by implementing an appropriate treatment process - Tr\_7) may be reused on App\_2. To sum up, the design decisions are: a) Connection of WWS\_1, b) alternate connections to App\_1, and c) treatment process from App\_1 to App\_2. The updated network is shown on Figure \ref{fig:case_3_to_be}.

%\begin{figure}[h!]    
%\centering
%\includegraphics[width=0.7\textwidth]{case_3_network_to_be.png}
%\caption{Updated network status of Chemical Industry A}
%\label{fig:case_3_to_be}
%\end{figure}

%For the operational phase, in order to validate the capability of our system, we will implement a set of experiments on the corresponding network status, derived from the design phase, under different possible values of input parameters, i.e. quality parameters and water flows, and for different values of the discretisation parameter $|K|$.

%Finally, we aim at comparing the two states of the network (before and after the implementation of the DSS) on a set of KPIs. The KPIs considered were designated by the case study and refer to the fresh water intake and the wastewater reuse. 

%Similarly to the Oil Refinery case, we work on both design and operational phases for Chemical Industry A. As mentioned in Section \ref{Section:usescenarioChemicalA}, there are three design decisions to be examined: i) the link from WWS\_1 on one of the treatment processes Tr\_1, Tr\_2 or Tr\_3, ii) the link from Tr\_1, Tr\_2 or Tr\_3 to App\_1 and iii) the most appropriate treatment process (Tr\_7) to treat effluent water from App\_1 and direct it to App\_2.

%We include these different options on the \textbf{optionsCompared} and \textbf{conflicts} fields. The number of trials $N$ is set to 500 and data are generated using the procedure described by Algorithm \ref{alg:alg_design}. 

Similarly to the Oil Refinery case, the number of trials $N$ is set to 500 and data are generated as described in Algorithm \ref{alg:alg_design}. To visualise our results, design decisions 2 and 3 are grouped (Table \ref{tab:chemA_setting_to_scenario}). Tr\_7v1 refers to the option with the greater conservation of flowrate but with lower quality improvement capability, while the opposite holds for Tr\_7v2. Table \ref{tab:chemA-res-design} presents the frequency of each option on the optimal solutions returned. The average execution of each trial is approximately 18 seconds, including the preprocessing stages. 

\begin{table}[h!]
\centering
 \begin{tabular}{||c || c ||} 
 \hline
 Setting & Decision\\ [0.5ex] 
 \hline\hline
 Set A & Tr\_1  $\rightarrow$ App\_1 + Tr\_7v1  \\ 
  Set B & Tr\_2  $\rightarrow$ App\_1  + Tr\_7v1\\
   Set C & Tr\_3  $\rightarrow$ App\_1  + Tr\_7v1 \\ 
  Set D& Tr\_1  $\rightarrow$ App\_1  + Tr\_7v2\\
   Set E & Tr\_2  $\rightarrow$ App\_1   + Tr\_7v2\\ 
  Set F& Tr\_3  $\rightarrow$ App\_1  + Tr\_7v2\\
 \hline
 \end{tabular}
 \caption{Grouping of Design Desisions for Chemical Industry A}
  \label{tab:chemA_setting_to_scenario}
\end{table}

\begin{table}[h!]
\centering
 \begin{tabular}{||c || c c c c c c| c||} 
 \hline
  & Set A & Set B &   Set C & Set D  & Set E & Set F & Sum\\ [0.5ex] 
 \hline\hline
 WWS\_1 $\rightarrow$ Tr\_1 & 61 &  72 & 39 & 3 & 9 &1 & 185\\ 
  WWS\_1 $\rightarrow$ Tr\_2 & 120 &  45 & 35 & 10 & 2 & 3 & 215\\
    WWS\_1 $\rightarrow$ Tr\_3 & 10 &  31 & 41 & 6 & 10 & 2 & 100 \\
    \hline
    Sum & 191 &  148 & 115 & 19 & 21 & 6 & 500\\    
 \hline
 \end{tabular}
 \caption{Design phase results for Chemical Industry A}
  \label{tab:chemA-res-design}
\end{table}

From Table \ref{tab:chemA-res-design}, it seems that the combination of Setting A (i.e. the link from Tr\_1 to App\_1 along with Tr\_7 v1) and the link from WWS\_1 to Tr\_2 is the most frequent choice of the optimisation algorithm, occurring on 120 trials (24\%). On the comparison of Tr\_7v1 (Settings A, B and C) and Tr\_v2 (Settings D, E and F), it is obvious that the former is preferred as it is utilised on a total of 454 (90.8 \%) versus 46 (9.2 \%) of the latter. This could be attributed to the higher conservation rate of Tr\_7v1, which allows for a greater amount of flowrate to be directed to App\_2. As a result, this leads to a greater percentage of the total flowrate demanded on App\_2, being recovered by the effluent of App\_1. In terms of the candidate links to App\_1, it seems that connecting Tr\_1 to App\_1 (Settings A and D) is more preferable with a total of 210 occurrences (42 \%) against the 169 (33.8 \%) of the second most frequent (Settings B and E). This shows that the quality parameters on the outlet of Tr\_1 are within the specified requirements of App\_1 and thus, it is possible to skip treatment process Tr\_4. Finally, regarding  WWS\_1, it is recommended to connect this wastewater stream with Tr\_2 as it is both identified on the most dominant combination with Setting A and on the total number of occurrences (215). Based on the above discussion, it seems that the most appropriate design decision is the combination of WWS\_1 $\rightarrow$ Tr\_2, the link of Tr\_2 to App\_1 and the employment of Tr\_7v1 between App\_1 and App\_2.

For the operational phase results, the above design decisions are employed, while adopting a similar procedure for the algorithmic evaluation as the one for the Oil Refinery (Table \ref{tab:chemA-res-oper}).

%By setting different values for $K$, we test the scalability of our method on datasets in terms of computational time, number of constraints and derived variables. We, again, consider the same set of representative values for $K\in \{100, 200, 500, 1000\}$. We generate a single dataset of 500 independent input instances and test the different values of $K$ on datasets with the same parameters. Computational results are shown in Table \ref{tab:chemA-res-oper}. Columns Vars, BVars and CVars correspond to the number of variables per type (B, C refer to Binary, Continuous, respectively) generated for each different value of $|K|$. Constraints refers to the number of constraints needed, Time (s) is the computational time in seconds and Gap refers to the optimality gap. 

\begin{table}[h!]
\centering
 \begin{tabular}{||c || c c c c c c||} 
 \hline
  $K$ & Vars & BVars &  CVars & Constraints & Time (s) & Gap (\%)\\ [0.5ex] 
 \hline\hline
 100 & 1506 & 1430 & 76 & 9293 & 2.58 & $<$ 0.7 \\ 
 200 &  2906 & 2830 & 76  &  18393 & 4.83 & $<$ 0.7 \\
  500 & 7106 & 7030 & 76 &  45693 & 17.83  & $<$ 0.6 \\ 
 1000 & 14106 & 14030 & 76 & 91193 & 32.66 & $<$ 0.5 \\
 \hline
 \end{tabular}
 \caption{Operational results on Chemical Industry A}
  \label{tab:chemA-res-oper}
\end{table}

As depicted on Table \ref{tab:chemA-res-oper}, optimal solutions for the operational phase are obtained in a few seconds. The time needed varies from approximately 2.5 s., when $K$ is set to 100 up to approximately 33 s., when $K$ is set to 1000. Similarly to our findings on the Oil Refinery case, we observe that increasing the value of $K$, affects only the number of binary variables (BVars) and Constraints, while the number of continuous variables (CVars) remains constant. Nevertheless, it is obvious that regardless of the value of $K$, our optimisation model remains highly efficient, producing optimal solutions within a few seconds. 

The final step on Chemical Industry A is to perform the comparison experiments, aiming to evaluate the enhancements on a set of critical KPIs proposed by the stakeholders of the plant. The compared networks correspond to the ones presented in Figures \ref{fig:case_3_as_is} (current network) and \ref{fig:case_3_to_be} (updated network). These KPIs correspond to i) the total \textbf{fresh water intake} used on the network, ii) the \textbf{reused} water and, iii) the percentage of \textbf{water losses}. 

\begin{table}[h!]
\centering
 \begin{tabular}{||c || c c c c ||} 
 \hline
  $Network$ & Total Int. & FW Int. & Reused \% & Losses \% \\ [0.5ex] 
 \hline\hline
 Current & 1782.125 & 1782.125 & 0 &  19.25 \\ 
 Updated & 1493.7 & 1468.4 & 10.58  & 9.01 \\
 \hline
 \end{tabular}
 \caption{Comparison of results on Chemical Industry A}
  \label{tab:chemA-res-comp}
\end{table}

The results of our comparison experiments are shown on Table \ref{tab:chemA-res-comp},  revealing the effectiveness of our approach from a practical perspective. We first highlight the obvious reduction in terms of the fresh water intake by 17.6 \%. Note that this percentage is close to the goal set in Section \ref{Section:CaseStudies} for that case study (20\%). The inclusion of WWS\_1, as well as the link of the effluent from App\_1 to App\_2 increases the reused flowrate, which contributes up to 10.58 \% to the final water used on the updated network. Finally, there is a significant reduce on the water losses throughout the network, originating from the decisions to bypass treatment processes that do not seem necessary for quality improvement (on certain instances). 

The above results highlight the merits of our DSS on Chemical Industry A. These benefits are established on: i) the design phase, where our method is employed in order to consult stakeholders on design decisions, ii) the operational phase, where our experiments have demonstrated the effectiveness (both in terms of computational time and solution quality) of our model and, iii) the comparison results, where it is numerically shown that our approach yields important enhancements on the KPIs proposed by the plant stakeholders.

\subsection{Chemical industry B} \label{Section:impactChemicalB}

The needs of this case study are limited on \textbf{design phase} decisions. The two main decisions to be made are: i) which outside treatment option is the less costly - yet sustainable for the plant and ii) up to which degree should the plant treat its wastewater stream before discharging into the respective agent. As a result, the DSS will be employed aiming at consulting plant stakeholders to invest upon the network design, with the goal of minimising the cost (i.e. tariffs imposed by the water providing actors) of the whole network operation.

In more detail, the network consists of one wastewater stream that needs to be treated, two applications where the purified water will be reused, and three different treatment options. Before the wastewater stream enters the network, there are three different possible treatment capabilities (Tr\_0); i) \textbf{low level treatment - LLT}, resulting in higher values of quality parameters, yet lower cost (due to the reduced use of chemicals), ii) \textbf{medium level treatment - MLT}, resulting in improved water quality, but greater costs and, iii) \textbf{high level treatment - HLT}, where water quality is enhanced at the greater degree possible, yet higher costs. In terms of the applications of the network, App\_1 demands a flowrate of greatly enhanced quality, while App\_2 is not as demanding regarding quality, but highly demanding in terms of flowrate - approximately 30 times greater than the one of App\_1.

\begin{figure}[h!]    
\centering
\includegraphics[scale=.5]{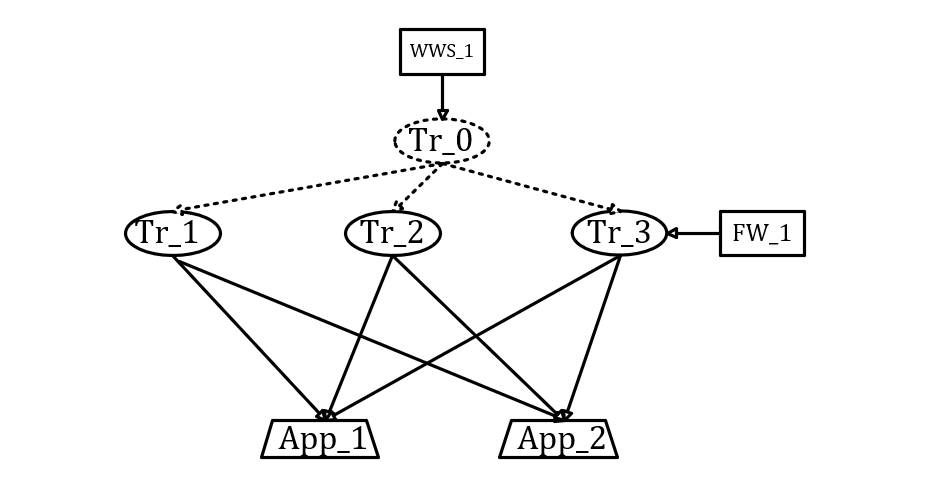}
\caption{Design network of Chemical Industry B}
\label{fig:case_2}
\end{figure}

Lastly, there are three possible treatment options, involving two external actors, the municipal wastewater treatment (Actor A) and a water purification system (Actor B). Option A (denoted by Tr\_1 on Figure \ref{fig:case_2}) consists of discharging wastewater to Actor A and then, receiving a predefined flow rate and water quality, to be directed to App\_1 and App\_2. Option A uses a fixed pricing scheme by inflicting the same cost regardless the quantity to be returned. Option B (Tr\_2) involves discharging wastewater to Actor B by offering a diverse pricing scheme, where the client is charged based on the water quality and flowrate. As a result, greater amount of flowrate corresponds to a reduced price (per $m^3$/h) - this price may be, also, reduced when the demanded flowrate is connected with lower quality requirements. Option C (Tr\_3) refers to the mixing of the wastewater stream with fresh water from Actor A (FW\_1) in the plant premises. On the one hand, the industry is not burdened with the discharging costs (though, there is an assumed investment cost), yet the requirements to enable this option are greater compared to Options A and B. Table \ref{tab:chemB_options_comp} presents a complete comparison of the candidate Options, summing up the above discussion.

\begin{table}[h!]
\centering
 \begin{tabular}{||c || c c c ||} 
 \hline
  Metric & Option A & Option B &  Option C\\ [0.5ex] 
 \hline\hline
 Discharge fee & Medium & High & Low  \\ 
 Requirements & Low & Medium & High \\
 Price to App\_1 & Low & High & Low  \\ 
 Price to App\_2 & Medium & Low & Medium \\
 \hline
 \end{tabular}
 \caption{Comparison of Options for Chemical Industry B}
  \label{tab:chemB_options_comp}
\end{table}

%Finally, note that from a strategic point of view, it is imperative that only a single option will be employed - thus, combination of different options is deemed invalid. 

Our approach consists of 2 different scenarios, aiming to determine the less costly combination of inbound wastewater and secondary treatment options. The first one is the \textbf{Baseline Scenario}, where the setting corresponds to the obligation of selecting a single secondary treatment Option. The second one, namely \textbf{All Options Available} is to determine whether the simultaneous inclusion of all options will yield different solutions. Table \ref{tab:chemB-res} presents the results, where each combination of options was used on the optimal solution.

%Each scenario is associated with 4 different experimental settings.

%\begin{itemize}
%    \item Setting A: Relax the quality requirements constraints with possible combinations among treatments i.e. there is the possibility to employ different options simultaneously
%    \item Setting B: Relax the quality requirements constraints, yet with the obligation of choosing only single options
%    \item Setting C: Enforce the quality requirements constraints allowing for possible combinations among treatments
%    \item Setting D: Enforce the quality requirements constraints with the obligation of choosing only single options.
%\end{itemize}

%Note that, for Settings A and B, there will be single solutions - as they are not dependent on the stochasticity of input quality parameters and costs are fixed per $m^3$/h. On the contrary, Settings C and D will be executed for $N$ iterations in order to determine the most robust option.  

\begin{table}[h!]
\centering
 \begin{tabular}{||c || c c c||} 
 \hline
  & Option A & Option B & Option C \\ [0.5ex] 
 \hline\hline
 LLT & 154 & 0 & 1 \\ 
 MLT & 27 & 0 & 66  \\
 HLT & 0 & 0 & 252  \\ 
 \hline
 \end{tabular}
 \caption{Design Phase Results for Chemical Industry B on the Baseline Scenario}
 \label{tab:chemB-res}
\end{table}

From Table \ref{tab:chemB-res}, it is obvious that Option B should not be considered, as it has not been deployed on any of the algorithm executions - possibly due to its greater discharge cost and the higher tariff on the water of enhanced quality. On the other hand, Option C is utlised on 319 trials (63.8 \%), with the majority of which being selected with the high level treatment (HLT). This is caused by the stricter necessary water quality requirements of Option C and as a result, the HLT is necessary. On the other hand, on some cases the mixing of wastewater with fresh water is not sufficient to satisfy the quality requirements of App\_1. Thus, it is preferable to mildly treat the water inside the plant (LLT or MLT) and then, opt for Option A (181 trials). As a result, it seems that for Chemical Industry B, the most appropriate combination of options would be the investment on internal plant infrastructure (Option C) with either medium (MLT) or high (HLT) level treatment. This diversification on treatment quality could be achieved by the addition of chemicals dynamically on the (MLT) in order to be transformed to (HLT), when quality requirements are not satisfied. 

Furthermore, it would be interesting to examine the `All Options Available' scenario, where the user may lift the `conflicts' constraints, which restrict the algorithm to select a single option. The corresponding results are shown on (Table \ref{tab:chemB-res-allopt}). 

\begin{table}[h!]
\centering
 \begin{tabular}{||c || c c c||} 
 \hline
  & Option A & Option B & Option C \\ [0.5ex] 
 \hline\hline
 LLT & 65 & 66 & 1 \\ 
 MLT & 17 & 74 & 62  \\
 HLT & 0 & 348 & 355 \\ 
 \hline
 \end{tabular}
 \caption{Design Phase Results for Chemical Industry B on the All Options Available Scenario}
  \label{tab:chemB-res-allopt}
\end{table}

First, note that the sum of the Table is greater than the number of trials (988 over 500) as the optimisation solutions for that scenario include - on the majority of the trials - a combination of two separate options. Numerical results change dramatically related to the ones presented on Table \ref{tab:chemB-res} in terms of the most preferable options. Most notably, Option B rises as a key option on this scenario being employed on 488 trials (97.6 \%). This is due to the alternative pricing scheme that Option B offers. More precisely, since the flowrates needed for App\_2 are much greater compared to the ones for App\_1, the price provided from Option B is the most competitive one. As a result, Option B is combined with either Option C - when the mixing of wastewater and fresh water returns eligible water quality to be directed to App\_1 - or with Option A on the opposite case.

The above results demonstrate how our DSS may be employed by an industrial user in order to assess both different investment possibilities and also, different pricing schemes, aiming at minimising the total operational cost. As a result, our tool provides the user the potential to examine different scenarios with the manipulation of pricing differentiations and even negotiate with the respective stakeholders.
\section{Concluding remarks} \label{Section:conlcusion}

We showcased the deployment of a Decision Support System in an industrial environment, evaluating its performance on network design and operational decision-making in an oil refinery and two distinct chemical industries. More precisely, after identifying the literature gap regarding DSS and optimisation algorithms in industrial settings, a generic data model was proposed, coupled by an efficient optimisation model. The proposed framework surpasses the limitations of current state-of-the-art literature, with its flexible design enhancing the method's applicability across varying networks, constraints, and objectives. This achievement can also be attributed to the more compact mathematical formulation, which offers a significant advantage over the problem-specific formulations commonly found in the literature. The aforementioned approach was applied to three case studies from diverse industries, each with distinct scopes and objectives, including chemical production and oil \& gas refining. As shown in Section \ref{Section:CaseStudies}, the first two aim at minimising either the fresh water intake (\textit{Chemical Industry A}) or the total cost (\textit{Chemical Industry B}), while the latter opts for maximising its wastewater treatment capability. The horizontal deployment of our methodology highlights its generic design, while our extensive experimentation setting, raises a number of non negligible benefits for the corresponding users. 

From a practical perspective, it is shown that all case studies were supported in making the most efficient design decisions aligned with their respective objectives. Especially for the \textit{Oil Refinery} and \textit{Chemical Industry A}, through the conducted experiments it is demonstrated that operational phase results are obtained within seconds, under different values of the discretisation number $K$. This automated procedure can significantly streamline the decision-making process for managers and plant owners, particularly during malfunction events where timely response is critical. In addition, the comparison with the current practice highlights the benefits of the new network structure with the critical freshwater intake objective reduced by 17.6\% for \textit{Chemical Industry A}, while the amount of wastewater reused in the \textit{Oil Refinery} increases to nearly 90\%. From a business perspective, deploying this system for real industrial needs involves three key roles: (a) IT personnel to configure and set up the system, (b) data engineers to integrate the required data from production systems into the solution, and (c) production experts (e.g., production/chemical engineers, line managers) to design and operate the water network.

On the other hand, it is important to mention that the required linearisation on our model exhibits few limitations. As already stated, there is a trade-off between the accuracy of results and computational time depending on the discretisation number $K$. Although in our work, the computational time burden imposed by higher values of $K$ is considered negligible (resulting in only a few extra seconds) for practical applications, it should be noted that in cases where even greater precision is required (e.g., setting $K$ to 10000), the problem may become intractable. Another concern may be that some real-world treatment processes operate under complicated equations, which hinder the explicit calculation of some of their corresponding attributes e.g. quality reductions rates ($RR_{ip}$). However, this challenge can be addressed through simulation and prediction techniques, which may generate the necessary data to serve as input for the DSS.

Concerning possible next steps, it seems interesting to focus on two complementary aspects: a) the design of a robust optimisation approach \cite{Bertsimas18, Bram15, Zhou23}, in order to take into consideration decision situations under uncertainty, identified for example on quality parameters and/or treatment efficiencies, and b) the application of our DSS on cases of networks where, apart from (waste-)water materials, such as by-products \cite{Cimren11}, raw materials \cite{Baxter16} or energy \cite{Badami19} are considered.

%This work gives floor to multiple possible expansions as the basis of new research to follow. First of all, an important task is to design appropriate methodologies, so as to tackle the most general variants of the WWMP problem. This rises as a severe challenge as it would extend the applicability of solvable flow networks by our optimization algorithm. On the other hand, dealing with uncertainty on industrial applications is a continuous dare for practitioners, since a common occurrence is the incapability of precisely determining the necessary network flow parameters. As a result, techniques of robust optimization (Ben-Tal, et al., 2008) should be taken under consideration in order to address possible uncertainties on e.g. incoming wastewater flow qualities or quantities. Finally, combining the optimization methods, as described above, with predictive analytics and real time monitoring could enable solving optimization problems not only on a limited time horizon, but, also, for multiple time cycles ahead.  
\appendix

\section{DSS User Interface} \label{sec:app_B}

\setcounter{figure}{0}

\begin{figure}[h]
\centering
\includegraphics[scale=.3]{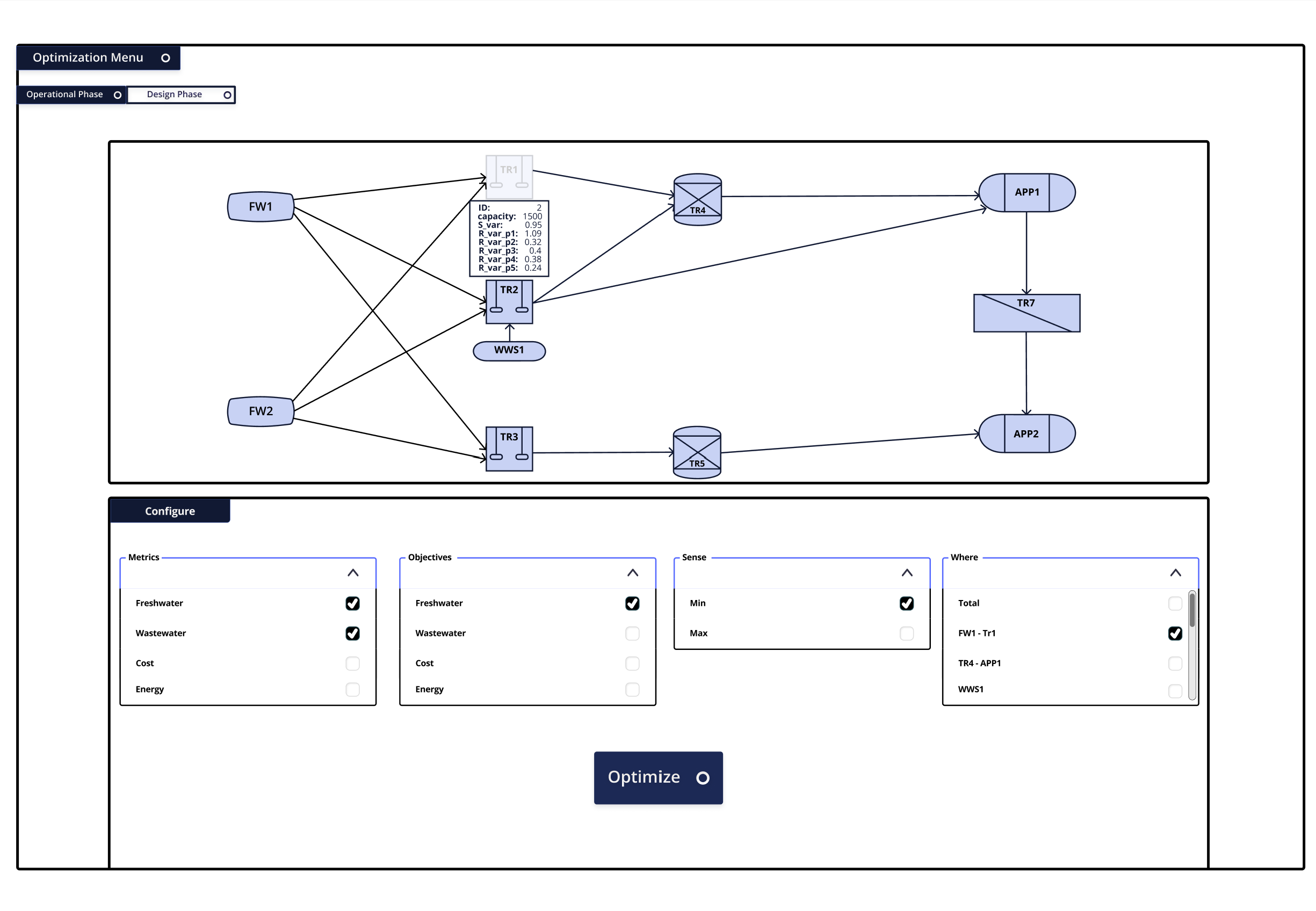}
\caption{User optimisation configuration under the operational phase}
\label{fig:oper_before}
\end{figure}

On \ref{sec:app_B}, we present indicative frames of the DSS User Interface. Due to space limitations, we focus on the \textit{Chemical Industry A} under the operational phase and on the \textit{Oil Refinery} under the design phase. However, these frames are applicable to all case studies across all phases, requiring only minor adjustments. Starting with the `operational phase', the user views the corresponding network and they are able to modify the attributes of each component by selecting it (e.g. for Tr\_1). Then, they proceed on configuring the optimisation model by selecting the corresponding choices on the dropdown menus: i) the metrics which are to be calculated from the model, ii) the optimisation objective, iii) the minimisation or maximisation of the objective (`sense') and, iv) the network components or edges whose operation should be optimised (`where'). The algorithm is triggered by clicking the `Optimize' button. The above discussion is depicted on Figure \ref{fig:oper_before}.

\begin{figure}[h]
\centering
\includegraphics[scale=.3]{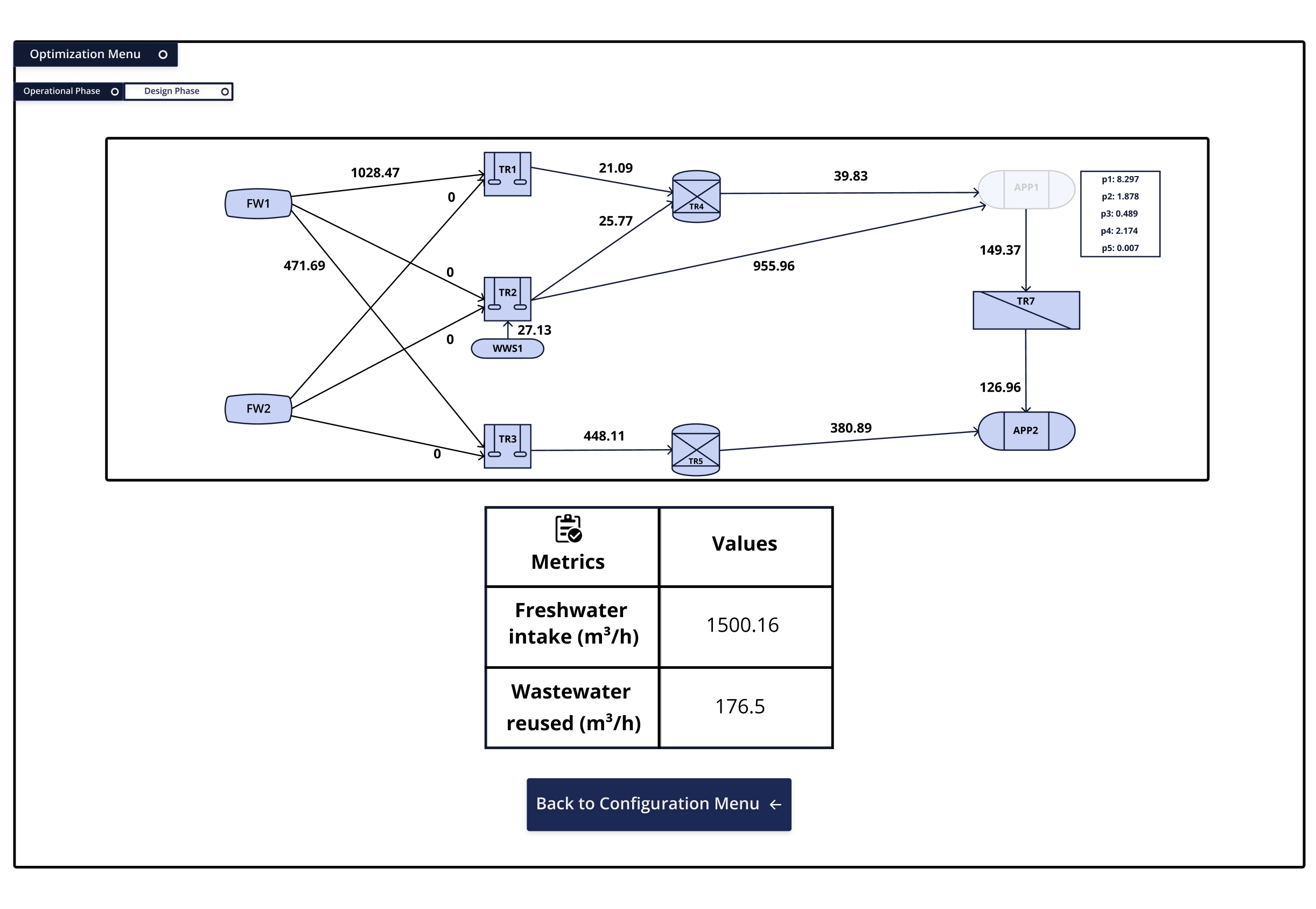}
\caption{User view of the results under the operational phase}
\label{fig:oper_after}
\end{figure}

Upon completion of the optimisation algorithm, the calculated flow rates for each edge are displayed. In addition to this feature, the user can view the calculated values for each network component (e.g., the inlet quality parameters on App1). The calculated metrics are also presented in a table format for further review. To configure a new optimisation request, the user can click the `Back to configuration menu' button. This process is illustrated in Figure \ref{fig:oper_after}.

In the design phase, we use the \textit{Oil Refinery} case study as an example. The user begins by configuring the network to be tested, adding new components and connecting them to existing ones using the toolbar on the left. 

\begin{figure}[h]
\centering
\includegraphics[scale=.3]{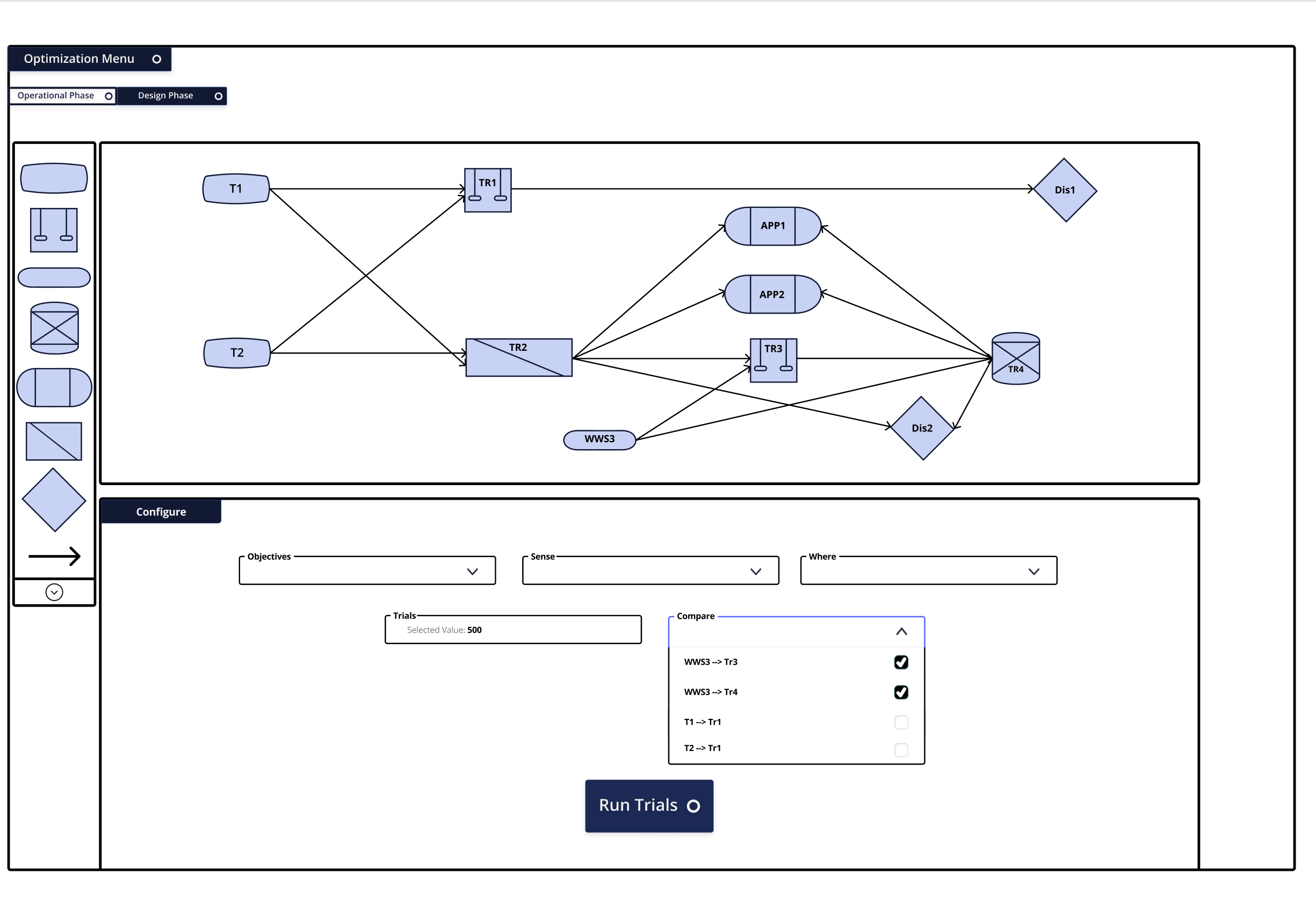}
\caption{User optimisation configuration under the design phase}
\label{fig:des_before}
\end{figure}

Additionally, the user can define sampling ranges and specify the distribution of data parameters for each relevant attribute. Once the network design is complete, the user configures the optimisation objective, selects the associated components (similar to the process shown in Figure \ref{fig:oper_before}) and specifies the number of trials to be executed. Finally, the user selects the conflicting options to be compared from the dropdown menu (`Compare'). The design optimisation is initiated by clicking the "Run Trials" button. This process is illustrated in Figure \ref{fig:des_before}.

\begin{figure}[h]
\centering
\includegraphics[scale=.3]{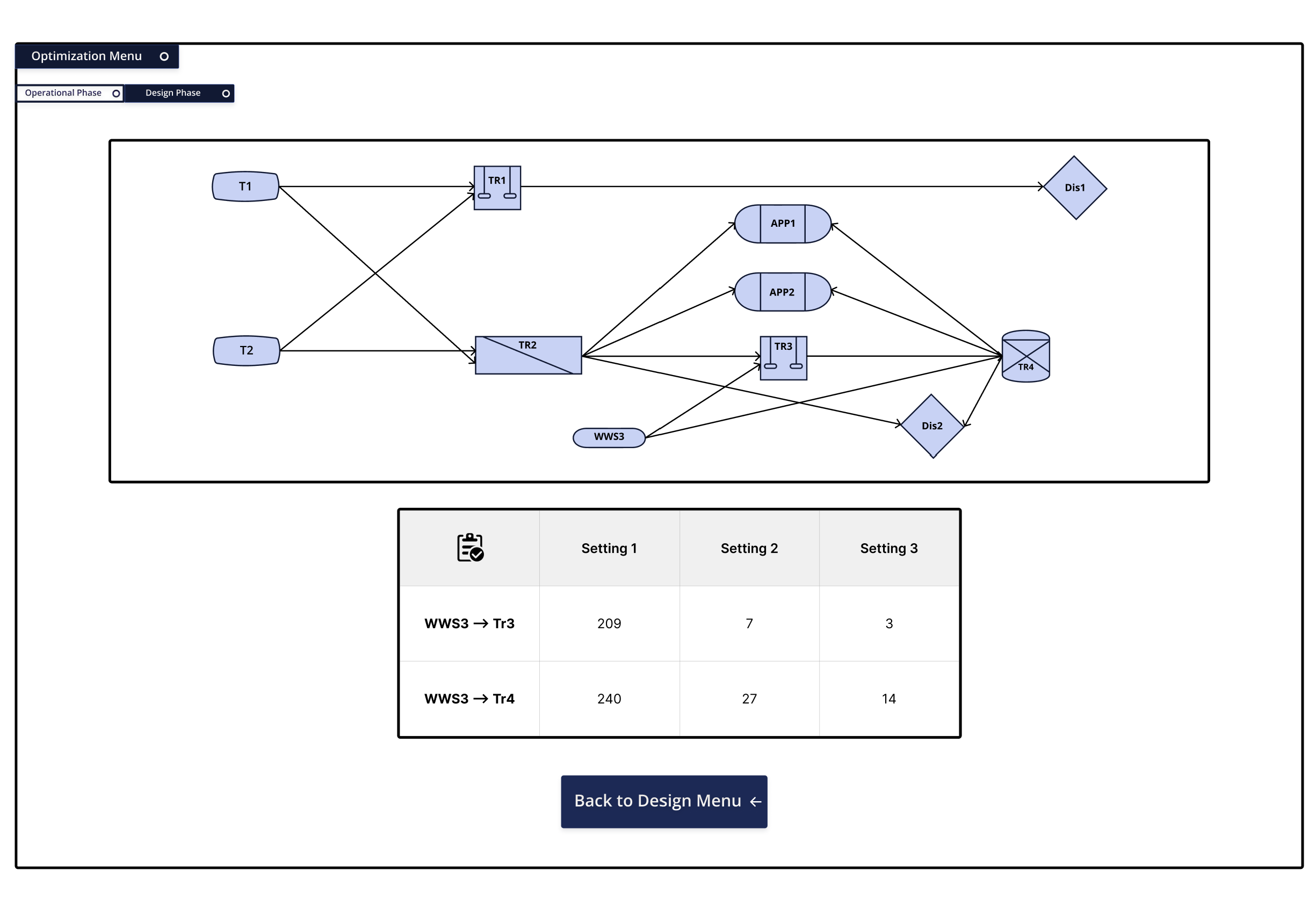}
\caption{User view of the results under the design phase}
\label{fig:des_after}
\end{figure}

Once the trials are completed, the user can view the frequency of each `compare' option appearing in the optimal solution across all trials. To configure a new design optimisation request, the user can click the `Back to Design Menu' button. An example frame is shown on Figure \ref{fig:des_after}.

\hfill \break
\textbf{Acknowledgments. }This work is supported by funds provided by the European Commission in the Horizon 2020 research and innovation programme AquaSPICE (Grant No. 958396) (\url{https://aquaspice.eu/}).

\setstretch{1.05}
\bibliographystyle{apa}

\newpage
\end{document}